\documentstyle[a4,12pt,amssymb,amstex,xypic]{amsart}
\xyoption{all}
\xyoption{arc}

\newtheorem{theo}{Theorem}[section]

\newtheorem{prop}[theo]{Proposition}

\newtheorem{lem}[theo]{Lemma}

\newtheorem{coro}[theo]{Corollary}

\def\et{{\mathrm{\acute{e}t}}}

\def\lalala{\calmath}

\def\lalala{\mathcal}

\def\FM{{\Bbb{F}}}

\def\QM{{\Bbb{Q}}}

\def\ZM{{\Bbb{Z}}}

\def\gGB{{\boldsymbol{\frak g}}}
\def\hGB{{\boldsymbol{\frak h}}}

\def\lGB{{\boldsymbol{\frak l}}}

\def\pGB{{\boldsymbol{\frak p}}}

\def\vGB{{\boldsymbol{\frak v}}}

\def\a{\alpha}
\def\b{\beta}

\def\g{\gamma}

\def\d{\delta}
\def\D{\Delta}

\def\ph{\varphi}
\def\ch{\chi}
\def\l{\lambda}

\def\s{\sigma}

\def\th{\theta}
\def\Th{\Theta}
\def\t{\tau}

\def\z{\zeta}

\def\Sigb{{\boldsymbol{\Sig}}}

\def\AC{{\lalala{A}}}

\def\DC{{\lalala{D}}}
\def\EC{{\lalala{E}}}
\def\FC{{\lalala{F}}}

\def\HC{{\lalala{H}}}

\def\KC{{\lalala{K}}}
\def\LC{{\lalala{L}}}
\def\MC{{\lalala{M}}}

\def\WC{{\lalala{W}}}
\def\XC{{\lalala{X}}}

\def\FCt{{\tilde{\lalala{F}}}}

\def\KCt{{\tilde{\lalala{K}}}}

\def\FCh{{\hat{\lalala{F}}}}

\def\NCB{{\boldsymbol{\lalala{N}}}}
\def\OCB{{\boldsymbol{\lalala{O}}}}

\def\SCB{{\boldsymbol{\lalala{S}}}}

\def\Ab{{\bold A}}
\def\Bb{{\bold B}}
\def\Cb{{\bold C}}

\def\Fb{{\bold F}}
\def\Gb{{\bold G}}
\def\Hb{{\bold H}}

\def\Lb{{\bold L}}
\def\Mb{{\bold M}}

\def\Ob{{\bold O}}
\def\Pb{{\bold P}}
\def\Qb{{\bold Q}}

\def\Sb{{\bold S}}
\def\Tb{{\bold T}}
\def\Ub{{\bold U}}
\def\Vb{{\bold V}}

\def\Xb{{\bold X}}
\def\Yb{{\bold Y}}
\def\Zb{{\bold Z}}

\def\pb{{\bold p}}

\def\zb{{\bold z}}

\def\fti{{\tilde{f}}}
\def\gti{{\tilde{g}}}

\def\uti{{\tilde{u}}}

\def\Kti{{\tilde{K}}}

\def\aha{{\hat{a}}}

\def\fha{{\hat{f}}}

\def\wdo{{\dot{w}}}

\def\fba{{\bar{f}}}

\def\Fbt{{\tilde{\Fb}}}

\def\Obt{{\tilde{\Ob}}}

\def\Xbt{{\tilde{\Xb}}}
\def\Ybt{{\tilde{\Yb}}}

\def\Xbh{{\hat{\Xb}}}
\def\Ybh{{\hat{\Yb}}}

\def\Sigb{{\boldsymbol{\Sigma}}}

\def\chit{{\tilde{\chi}}}

\def\sigt{{\tilde{\s}}}

\def\thet{{\tilde{\theta}}}

\def\alpba{{\bar{\alpha}}}
\def\betba{{\bar{\beta}}}

\def\piba{{\bar{\pi}}}

\def\ad{\mathop{\mathrm{ad}}\nolimits}
\def\Ad{\mathop{\mathrm{Ad}}\nolimits}

\def\diag{\mathop{\mathrm{diag}}\nolimits}

\def\End{\mathop{\mathrm{End}}\nolimits}

\def\Id{\mathop{\mathrm{Id}}\nolimits}

\def\Im{\mathop{\mathrm{Im}}\nolimits}

\def\Ind{\mathop{\mathrm{Ind}}\nolimits}

\def\INT{\mathop{\mathrm{int}}\nolimits}

\def\Ker{\mathop{\mathrm{Ker}}\nolimits}

\def\reg{{\mathrm{reg}}}

\def\Res{\mathop{\mathrm{Res}}\nolimits}

\def\Stab{\mathop{\mathrm{Stab}}\nolimits}

\def\supp{\mathop{\mathrm{supp}}\nolimits}

\def\Tr{\mathop{\mathrm{Tr}}\nolimits}

\def\uni{{\mathrm{uni}}}

\def\tete#1{\par\leavevmode\makebox[0.7cm]{$(\mathrm{#1})$}}

\def\to{\rightarrow}
\def\longto{\longrightarrow}
\def\injto{\hookrightarrow}

\def\longmapright#1{\hspace{0.3em}\smash{
     \mathop{\longrightarrow}\limits^{#1}}\hspace{0.3em}}

\def\fonction#1#2#3#4#5{\begin{array}{rccc}
{#1} : & {#2} & \longto & {#3} \\
& {#4} & \longmapsto & {#5} 
\end{array}}

\def\fonctio#1#2#3#4{\begin{array}{ccc}
{#1} & \longto & {#2} \\
{#3} & \longmapsto & {#4} 
\end{array}}

\def\ci{\circ}
\def\pr{\prime}

\def\incl{\hspace{0.05cm}{\subset}\hspace{0.05cm}}
\def\notincl{\hspace{0.05cm}{\not\subset}\hspace{0.05cm}}

\def\fq{\FM_q}

\def\ql{{\QM_\el}}
\def\qlb{{\overline{\QM}_\el}}

\def\DS{\displaystyle}
\def\SS{\scriptstyle}

\def\fin{~$\SS \blacksquare$}
\def\finl{~$\SS \square$}

\def\el{\ell}

\def\matrice#1{\left(\begin{array}{ccccccccccccccccccc}#1\end{array}\right)}

\def\lexp#1#2{\kern\scriptspace\vphantom{#2}^{#1}\kern-\scriptspace#2}
\def\le{\hspace{0.1em}\mathop{\leqslant}\nolimits\hspace{0.1em}}
\def\ge{\hspace{0.1em}\mathop{\geqslant}\nolimits\hspace{0.1em}}

\mathchardef\lllllll="3278
\def\SEC{$\lllllll$}
\mathchardef\inferieur="321E
\mathchardef\superieur="321F
\def\arobas{\char'100}

\def\example{\noindent{\sc{Example}~-} }

\def\sec{\section}
\def\sub{\subsection}

\def\bi{\bigskip}
\def\med{\medskip}
\def\sma{\smallskip}
\def\equanum{\begin{equation}\begin{array}{r@{=}l}}
\def\endequanum{\end{array}\end{equation}}
\def\eqna{\begin{eqnarray*}}
\def\endeqna{\end{eqnarray*}}

\def\mor{morphism }

\def\endo{endomorphism }

\def\iso{isomorphism }
\def\isos{isomorphisms }
\def\auto{automorphism }

\def\proof{\noindent{\sc{Proof}~-} }
\def\rem{\noindent{\sc{Remark}~-} }

\def\borel{Borel subgroup }

\def\para{parabolic subgroup }
\def\paras{parabolic subgroups }
\def\levi{Levi subgroup }
\def\levis{Levi subgroups }

\def\car{character }

\def\irr{irreducible }

\def\resp{respectively }

\def\cf{{\it cf.}~}

\def\equat{\refstepcounter{theo}$$~}
\def\endequat{\leqno{\boldsymbol{(\arabic{section}.\arabic{theo})}}~$$}

\newcounter{numero}[section]

\def\na{\nabla}

\def\ins{{\mathrm{ins}}}

\def\mini{{\mathrm{min}}}

\def\remark#1{{\refstepcounter{theo}\label{#1}\noindent\sc Remark  
\arabic{section}.\arabic{theo} - }}
\def\example#1{{\refstepcounter{theo}\label{#1}\noindent\sc Example 
\arabic{section}.\arabic{theo} - }}

\begin{document}

\begin{centerline}{\Large \bf Actions of relative Weyl groups I}\end{centerline}

\bi

\begin{centerline}{\sc C\'edric Bonnaf\'e}\end{centerline}

\bi

\begin{centerline}{\today}\end{centerline}

\bi

\bi

\begin{centerline}{\sc General introduction}\end{centerline}

\bi

A theorem of Digne, Lehrer and Michel says that Lusztig restriction of 
a Gel'fand-Graev \car of a finite reductive group $\Gb^F$ is still a Gel'fand-Graev 
character \cite[Theorem 3.7]{DLM2}. However, 
an ambiguity remains on the character obtained (whenever the center of $\Gb$ 
is not connected, there are several Gel'fand-Graev characters). 
The original aim of this series of papers was to drop this ambiguity. 
For this, we needed to study more deeply the structure of the 
endomorphism algebra of an induced cuspidal character sheaf~: 
for instance, we wanted to follow the action of a Frobenius 
endomorphism through this algebra.

This lead us to this first part, in which we develop another approach 
for computing explicitly this endomorphism algebra. One of the main goal 
is to construct another isomorphism between this endomorphism algebra and 
the group algebra of the 
relative Weyl group involved (a first one was constructed by 
Lusztig \cite[Theorem 9.2]{luicc}). By comparing both isomorphisms, 
we get some immediate consequences for finite reductive groups. 
Note that the results of this part are valid for any 
cuspidal local system supported by a unipotent class and have 
a chance to be useful for computing values of characters 
at unipotent elements. 

In the second forthcoming part, we will restrict our attention to the case where 
the cuspidal local system is supported by the regular unipotent 
class. We will then be able to compute explicitly the generalized 
Springer correspondence through this new isomorphism. As an application 
of these (sometimes fastidious) computations, we will get 
the desired more precise version of Digne, Lehrer and Michel's theorem. 
It must be said that this final result is valid only whenever the cardinality 
of the finite field is large enough.

\bi

\begin{centerline}{\sc Introduction to this first part}\end{centerline}

\bi

Let $\Gb$ be a connected reductive group defined over an algebraically 
closed field $\FM$, let $\Lb$ be a \levi of a \para of $\Gb$, let $\Cb$ 
be a unipotent class of $\Lb$, and let $v \in \Cb$. We first explain how the action of 
the finite group $W_\Gb(\Lb,\Cb)=N_\Gb(\Lb,\Cb)/\Lb$ on some varieties 
introduced by Lusztig \cite[\SEC 3 and 4]{luicc} can be extended 
``by density'' to some slightly bigger varieties (see \SEC \ref{sub relatif}). 
We then generalize this construction to extend the action of 
$W_\Gb(\Lb,v) =N_\Gb(\Lb,v)/C_\Lb^\ci(v)$ on other varieties 
covering the previous ones (see \SEC \ref{sub nor}). We show 
that, in most cases, the description of some stabilizers 
is entirely given by a morphism 
$\ph_{\Lb,v}^\Gb : W_\Gb(\Lb,\Cb) \to A_\Lb(v)=C_\Lb(v)/C_\Lb^\ci(v)$ 
(see \SEC \ref{02} and \ref{01}). We also provide some reduction 
arguments to compute explicitly the morphisms $\ph_{\Lb,v}^\Gb$ 
(see \SEC \ref{ele sec}).

From \SEC \ref{sec endo} to the end, we assume that $\Cb$ supports 
a cuspidal local system $\EC$ (we denote by $\z$ the 
\car of the finite group $A_\Lb(v)$ associated to $\EC$). Let 
$K$ denote the perverse sheaf obtained from the datum $(\Cb,\EC)$ 
by parabolic induction \cite[4.1.1]{luicc}, and let $\AC$ denote 
its endomorphism algebra. In this case, $W_\Gb(\Lb,\Cb)$ is 
equal to $W_\Gb(\Lb)$ (by \cite[Theorem 9.2]{luicc}) and is 
isomorphic to $W_\Gb^\ci(\Lb,v)=N_\Gb(\Lb) \cap C_\Gb^\ci(v)/C_\Lb^\ci(v)$. 
Lusztig \cite[Theorem 9.2]{luicc} 
constructed a canonical isomorphism 
$\Th : \qlb W_\Gb^\ci(\Lb,v) \to \AC$. The aim of \SEC \ref{sec endo} 
is to construct an other explicit isomorphism $\Th' : \qlb W_\Gb^\ci(\Lb,v) \to \AC$.
For this purpose, we use the varieties previously introduced in \SEC \ref{mor phi}. 
It turns out that $\Th$ and $\Th'$ differs by a linear character 
$\g_{\Lb,v,\z}^\Gb$ (see Corollary \ref{epsilon}). Moreover, 
whenever $\z$ is linear, and whenever some mild hypothesis 
are satisfied, we have $\g_{\Lb,v,\z}=\z \ci \ph_{\Lb,v}^\Gb$. 

In \SEC \ref{part finite}, we assume further that $\FM$ is an algebraic 
closure of a finite field and that $\Gb$ is endowed with a Frobenius 
endomorphism. We then explain what kind of precisions may be obtained 
by using the previous results about characteristic functions of 
character sheaves.

In the future part II, we will assume throughout that $v$ is regular. 
Under this hypothesis, 
we will compute explicitly the morphisms $\ph_{\Lb,v}^\Gb$ 
in sufficiently many cases to be able to apply successfully this 
work to Gel'fand-Graev characters. We then follow Digne, Lehrer and 
Michel's method~: from the knowledge of $\g_{\Lb,v,\z}$ and from 
the explicit nature of the isomorphism $\Th'$, we get the 
desired precision.

\bi

\def\isol{{\mathrm{iso}}}

\section*{Notation}~

\med

\noindent{\bf Fields, varieties, sheaves.} We fix an algebraically 
closed field $\FM$ and we denote by $p$ its 
characteristic. All algebraic varieties and all algebraic groups will be 
considered over $\FM$. We also fix a prime number $\el$ different from $p$. 
Let $\qlb$ denote an algebraic 
closure of the $\el$-adic field $\ql$. 

If $\Xb$ is an algebraic variety (over $\FM$), we also denote 
by $\qlb$ the constant $\el$-adic sheaf associated to $\qlb$ (if necessary, 
we denote it by $(\qlb)_\Xb$). By a constructible 
sheaf (respectively a local system) on $\Xb$ we mean 
a constructible $\qlb$-sheaf (respectively a $\qlb$-local system).
Let $\DC\Xb$ denote the 
bounded derived category of constructible sheaves on $\Xb$. If $K \in \DC\Xb$ 
and $i \in \ZM$, we denote by $\HC^iK$ the $i$-th cohomology sheaf of $K$ 
and if $x \in \Xb$, then $\HC^i_x K$ denotes the stalk at $x$ of the 
constructible sheaf $\HC^iK$. If $K \in \DC \Xb$, 
we denote by $DK$ its Verdier dual. If $\LC$ is a constructible sheaf on $\Xb$, 
we identify it with its image in $\DC\Xb$, that is 
the complex concentrated in degree $0$ whose $0$th term is $\LC$. 

Let $K \in \DC \Xb$. We say that $K$ is a {\it perverse sheaf} if the following two 
conditions hold~:
$${\mathrm{(a)}}\hspace{2.3cm}\forall i \in \ZM,~\dim \supp \HC^i K \le -i,$$
$${\mathrm{(b)}}\hspace{2cm}\forall i \in \ZM,~\dim \supp \HC^i D K \le -i.$$
We denote by $\MC\Xb$ the full subcategory of $\DC\Xb$ whose objects are perverse 
sheaves~: this is an abelian category \cite[2.14, 1.3.6]{BBD}.

Let $\Yb$ be a locally closed, smooth, irreducible subvariety of $\Xb$ and let 
$\LC$ be a local system on $\Yb$. We denote by $IC(\overline{\Yb},\LC)$ the 
intersection cohomology complex of Deligne-Goresky-MacPherson of $\overline{\Yb}$ 
with coefficients in $\LC$. We often identify $IC(\overline{\Yb},\LC)$ with 
its extension by zero to $\Xb$~; $IC(\overline{\Yb},\LC)[\dim \Yb]$ is a perverse sheaf 
on $\Xb$. 

\bi

\noindent{\bf Algebraic groups.} If $\Hb$ is a linear algebraic group, we will denote 
by $\Hb^\ci$ the neutral component of $\Hb$, by $\Hb_\uni$ the closed 
subvariety of $\Hb$ consisting of unipotent elements of $\Hb$, and 
by $\Zb(\Hb)$ the center of $\Hb$. If $h \in \Hb$, then 
$A_\Hb(h)$ denotes the finite group $C_\Hb(h)/C_\Hb^\ci(h)$, $(h)_\Hb$ denotes the 
conjugacy class of $h$ in $\Hb$, and $h_s$ (\resp $h_u$) denotes 
the semisimple (\resp unipotent) part of $h$. If $\hGB$ is 
the Lie algebra of $\Hb$, we denote by $\Ad h : \hGB \to \hGB$ the 
differential at $1$ of the \auto $\Hb \to \Hb$, $x \mapsto \lexp{h}{x}=hxh^{-1}$.

If $\Xb$ and $\Yb$ are varieties, and if $\Xb$ (\resp $\Yb$) is 
endowed with an action of $\Hb$ on the right (\resp left), 
then we denote, when it exists, $\Xb \times_\Hb \Yb$ the quotient 
of $\Xb \times \Yb$ by the diagonal left action of $\Hb$ given by 
$h.(x,y)=(xh^{-1},hy)$ for any $h \in \Hb$ and $(x,y) \in \Xb \times \Yb$. 
If $(x,y) \in \Xb \times \Yb$, and if $\Xb \times_\Hb \Yb$ exists, 
we denote by $x *_\Hb y$ the image of $(x,y)$ in $\Xb \times_\Hb \Yb$ 
by the canonical morphism.

Finally, if $X_1$,\dots, $X_n$ are subsets or elements of $\Hb$, we denote 
by $N_\Hb(X_1,\dots,X_n)$ the intersection of the normalizers $N_\Hb(X_i)$ 
of $X_i$ in $\Hb$ ($1 \le i \le n$). 

\bi

\noindent{\bf Reductive group.} We fix once and for all a connected reductive algebraic 
group $\Gb$. We fix a \borel $\Bb$ of $\Gb$ and a maximal torus $\Tb$ of $\Bb$. 
We denote by $X(\Tb)$ (\resp $Y(\Tb)$) the lattice of rational characters 
(\resp of one-parameter subgroups) of $\Tb$. 
Let $W=N_\Gb(\Tb)/\Tb$. 
Let $\Phi$ denote the root system of $\Gb$ relative to $\Tb$ and let 
$\Phi^+$ (\resp $\D$) denote the set of positive roots 
(\resp the basis) of $\Phi$ associated to $\Bb$. 
For each root $\a \in \Phi$, we denote by 
$\Ub_\a$ the one-parameter unipotent subgroup of $\Gb$ 
normalized by $\Tb$ associated to $\a$. 

We also fix in this paper 
a \para $\Pb$ of $\Gb$ and a \levi $\Lb$ of $\Pb$. We denote by $\pi_\Lb : \Pb \to \Lb$ 
the canonical projection with kernel $\Vb$, the unipotent radical of $\Pb$. 
We denote by $\Phi_\Lb$ the root system 
of $\Lb$ relative to $\Tb$~; we have $\Phi_\Lb \incl \Phi$. Finally, $W_\Lb$ 
denotes the Weyl group of $\Lb$ relative to $\Tb$. 

\bi

\sec{Preliminaries\label{sec steinberg}}~

\med

\sub{Centralizers} We start this subsection by recalling two well-known 
results on centralizers of elements in reductive groups. The first one is 
due to Lusztig \cite[Proposition 1.2]{luicc}, while the second one is due 
to Spaltenstein \cite[Proposition 3]{holt}.

\bi

\begin{lem}[{\bf Lusztig}]\label{dim p}
$(1)$ Let $l \in \Lb$ and $g \in \Gb$. Then 
$$\dim \{x \Pb~|~x^{-1}gx \in (l)_\Lb .\Vb\} \le 
{1 \over 2}(\dim C_\Gb(g)-\dim C_\Lb(l)).$$

\tete{2} If $g \in \Pb$, then $\dim C_\Pb(g) \ge \dim C_\Lb(\pi_\Lb(g))$.
\end{lem}

\bi

\begin{lem}[{\bf Spaltenstein}]\label{spa connexe}
If $l \in \Lb$, then $C_\Vb(l)$ is connected.
\end{lem}

\bi

We will now give several applications of the two previous lemmas. 
We first need the following technical result~:

\bi

\begin{lem}\label{equiv class}
Let $l \in \Lb$. Then the following are equivalent~:

\tete{a} $C_\Gb^\ci(l) \incl \Lb$~;

\tete{b} $C_\Gb^\ci(l_s) \incl \Lb$~;

\tete{c} $C_\Vb(l)=\{1\}$~;

\tete{d} $\dim C_\Vb(l)=0$.
\end{lem}

\bi

\proof It is clear that (b) implies (a), and that (a) implies (d). 
Moreover, by Lemma \ref{spa connexe}, (c) is equivalent to (d). 
It remains to prove that (c) implies (b). 

For this, let $s$ (\resp $u$) denote the semisimple (\resp unipotent) 
part of $l$, and assume that $C_\Gb^\ci(s) \notincl \Lb$. We want 
to prove that $C_\Vb(l) \not= \{1\}$. Without loss of generality, 
we may, and we will, assume that $s \in \Tb$ and $u \in \Ub \cap \Lb$. 

Let $\Gb'=C_\Gb^\ci(s)$, $\Ub'=C_\Ub(s)$, $\Vb'=C_\Vb(s)$ and $\Bb'=C_\Bb(s)$. 
Then, by \cite[Corollary 11.12]{borel}, $u \in \Gb'$. 
Moreover, by Lemma \ref{spa connexe}, $\Ub'$ and $\Vb'$ are connected. 
So $\Bb'=\Tb.\Ub'$ is connected. Let $\Phi_s$ denote the root system 
of $\Gb'$ relative to $\Tb$, and let $\Phi_s^+$ 
be the positive root system associated to 
the Borel subgroup $\Bb'$ of $\Gb'$. Let $\a_s$ denote the highest 
root of $\Phi_s$ with respect to $\Phi_s^+$. 
Since $\Gb' \notincl \Lb$, we have $\a_s \not\in \Phi_\Lb$. 
But $\Ub_{\a_s}$ is central in $\Ub'$, so 
$\Ub_{\a_s} \incl C_{\Vb'}(u)=C_\Vb(l)$. 
Therefore, $C_\Vb(l) \not= \{1\}$. The proof of Lemma 
\ref{equiv class} is now complete.\fin

\bi

Now, let $\OCB$ be the set of elements $l \in \Lb$ such that 
$C_\Vb(l)=\{1\}$. 

\bi

\begin{lem}\label{U conjugue}
The set $\OCB$ is a dense open subset of $\Lb$ and the map
$$\fonctio{\OCB \times \Vb}{\OCB.\Vb}{(l,x)}{xlx^{-1}}$$
is an isomorphism of varieties.
\end{lem}

\bi

\proof The group $\Vb$ acts on $\Gb$ by conjugation. So, by 
\cite[Proposition 1.4]{humphreys}, the set 
$$\NCB=\{g \in \Gb~|~\dim C_\Vb(g) =0 \}$$
is an open subset of $\Gb$. Therefore, $\NCB \cap \Lb$ is an open 
subset of $\Lb$. Moreover, $\NCB \cap \Lb$ is not empty 
since any $\Gb$-regular element of $\Tb$ belongs 
to $\NCB \cap \Lb$. But, by Lemma \ref{spa connexe}, $\NCB \cap \Lb=\OCB$. 
This proves the first assertion of the lemma.

Now, let $f$ denote the \mor defined in Lemma \ref{U conjugue}. 
Let $l \in \OCB$. Then the map $f_l : \Vb \to l.\Vb$, $x \mapsto xlx^{-1}$ 
is injective by definition of $\OCB$, and its image is closed because it 
is an orbit under a unipotent group \cite[Proposition 4.10]{borel}. 
By comparing dimensions, we get that $f_l$ is bijective. 
As this holds for every $l \in \OCB$, $f$ is bijective.

Moreover, the variety $\OCB.\Vb \simeq \OCB \times \Vb$ is smooth. 
Hence, to prove that $f$ is an isomorphism, we only need to prove, 
by \cite[Theorems AG.17.3 and AG.18.2]{borel}, 
that the differential $(df)_{(l,1)}$ is surjective for 
some $l \in \OCB$. 

Now, let $t \in \Tb$ be a $\Gb$-regular element (so that $t \in \OCB$). 
The tangent space to $\OCB$ at $t$ may be identified with the Lie algebra 
$\lGB$ of $\Lb$ via the translation by $t$. By writing 
$f(l,x)=l.(l^{-1}xlx^{-1})$ for every $(l,x) \in \OCB \times \Vb$, 
the differential $(df)_{(t,1)}$ may be identified with the map 
$$\fonctio{\lGB \oplus \vGB}{\lGB \oplus \vGB}{l \oplus x}{l \oplus 
(\Ad t^{-1} - \Id_\vGB)(x).}$$
Here, $\vGB$ denotes the Lie algebra of $\Vb$. 
The bijectivity of $(df)_{(t,1)}$ follows immediately from the fact 
that the eigenvalues of $\Ad t^{-1}$ are equal to $\a(t)^{-1}$ for 
$\a \in \Phi^+-\Phi_\Lb$, so they are different from $1$ by 
the regularity of $t$.\fin

\bi

Lemma \ref{U conjugue} implies immediately the following result~:

\bi

\begin{coro}\label{U conjugue bis}
Let $\SCB$ be a locally closed subvariety of $\OCB$. Then the map
$$\fonctio{\SCB \times \Vb}{\SCB.\Vb}{(l,x)}{xlx^{-1}}$$
is an \iso of varieties.
\end{coro}

\bi

\noindent{\sc Notation - } If $\SCB$ is a locally closed subvariety 
of $\Lb$, we denote by $\SCB_\reg$ (or $\SCB_{\reg,\Gb}$ if there is 
some ambiguity) the open subset $\SCB \cap \OCB$ of $\SCB$. It might 
be empty.\finl

\bi

\sub{Steinberg map} Let $\na : \Gb \to \Tb/W$ be the 
Steinberg map. Let us recall its definition.
If $g \in \Gb$, then $\na(g)$ is the intersection of $\Tb$ with the conjugacy class 
of the semisimple part of $g$. Then $\na$ is a \mor of varieties 
\cite[\SEC{6}]{steinreg}. To compute the Steinberg map, we need to determine 
semisimple parts of elements of $\Gb$. In our situation, the following well-known 
lemma will be useful \cite[5.1]{luicc}~:

\bi

\begin{lem}\label{V conjugue}
If $g \in \Pb$, then the semisimple part of $g$ is 
$\Vb$-conjugate to the semisimple part of $\pi_\Lb(g)$. In particular, 
$\na(g)=\na(\pi_\Lb(g))$.
\end{lem}

\bi

\proof Let $s$ be the semisimple part of $g$. Then the semisimple part of 
$\pi_\Lb(g)$ is $\pi_\Lb(s)$. But, $s$ belongs to some \levi $\Lb_0$ of 
$\Pb$. Let $x \in \Vb$ be such that $\Lb=\lexp{x}{\Lb_0}$. 
Then $xsx^{-1} \in \Lb$ and $xsx^{-1}$ is the semisimple part of 
$xgx^{-1}$. Therefore, the semisimple part of 
$\pi_\Lb(xgx^{-1})=\pi_\Lb(g)$ is $\pi_\Lb(xsx^{-1})=xsx^{-1}$.\fin

\bi

%
%
%

\bi

\begin{lem}\label{nabla -2}
$\na(\Zb(\Lb)^\ci)$ is a closed subset of $\Tb/W$ and $\na(\Zb(\Lb)^\ci_\reg)$ 
is an open subset of $\na(\Zb(\Lb)^\ci)$. 
\end{lem}

\bi

\proof The restriction of $\na$ to $\Tb$ is a finite quotient morphism. 
In particular, it is open and closed. Since it is closed, $\na(\Zb(\Lb)^\ci)$ 
is a closed subset of $\Tb/W$. Since it is open, 
$\na(\Tb_\reg)$ is an open subset of $\Tb/W$. 
But $\na(\Zb(\Lb)^\ci_\reg)=\na(\Tb_\reg) \cap \na(\Zb(\Lb)^\ci)$. 
So the Lemma \ref{nabla -2} is proved.\fin

\bi 

Let $\na_\Lb : \Lb \to \Tb/W_\Lb$ denote the Steinberg map for the 
group $\Lb$. By Lemma \ref{equiv class}, we have
\equat\label{O T}
\OCB=\na_\Lb^{-1}(\Tb_\reg/W_\Lb).
\endequat

\bi

\sub{A family of morphisms\label{sub fam}} If $\SCB$ is a locally closed subvariety 
of $\Lb$ stable under conjugation by $\Lb$, then $\SCB.\Vb$ is a 
locally closed subvariety of $\Pb$ stable under conjugation by $\Pb$. 
We can therefore consider the quotients $\Gb \times_\Lb \SCB$ and 
$\Gb \times_\Pb \SCB.\Vb$. With these varieties are associated the 
maps $\Gb \times_\Lb \SCB \to \Gb$, $g *_\Lb l \mapsto glg^{-1}$ 
and $\Gb \times_\Pb \SCB.\Vb \to \Gb$, $g *_\Pb x \mapsto gxg^{-1}$~: 
they are well-defined morphisms of varieties. 

\bi

\remark{identification} If $\SCB$ is contained in $\OCB$, then 
the map $\Gb \times_\Lb \SCB \to \Gb \times_\Pb \SCB.\Vb$, 
$g *_\Lb l \mapsto g *_\Pb l$ is an isomorphism of varieties 
(by Corollary \ref{U conjugue bis}).\finl

\bi

The next result is well-known~:

\bi

\begin{lem}\label{proj mor}
The map $\Gb \times_\Pb \Pb \to \Gb$, $g *_\Pb x \mapsto gxg^{-1}$ is a projective 
surjective morphism of varieties. In particular, if $\Fb$ is a closed subvariety 
of $\Pb$ stable under conjugation by $\Pb$, then the map 
$\Gb *_\Pb \Fb \to \Gb$, $g *_\Pb x \mapsto gxg^{-1}$ is a projective 
morphism.
\end{lem}

\bi

\proof Let $\Xbt=\{(x,g\Pb) \in \Gb \times \Gb/\Pb~|~g^{-1}xg \in \Pb\}$. 
Then $\Xbt$ is a closed subvariety of $\Gb \times \Gb/\Pb$. Moreover, 
the variety $\Gb/\Pb$ is projective. Therefore, the projection 
$\pi : \Xbt \to \Gb$, $(x,g\Pb) \mapsto x$ is a projective morphism. 
Since every element of $\Gb$ is conjugate to an element of $\Bb$, 
$\pi$ is surjective.

But the maps $\Gb \times_\Pb \Pb \to \Xbt$, $g *_\Pb x \mapsto (gxg^{-1},g\Pb)$ 
and $\Xbt \to \Gb \times_\Pb \Pb$, $(x,g\Pb) \mapsto g *_\Pb g^{-1}xg$ 
are morphisms of varieties which are inverse of each other. 
Moreover, through these isomorphisms, the map constructed in Lemma \ref{proj mor} 
may be identified with $\pi$. The proof is now complete.\fin

\bi

The next result might be known but the author have never 
seen such a statement.

\bi

\begin{lem}\label{lissable}
The morphisms of varieties 
$$\fonctio{\Gb \times_\Lb \OCB}{\Gb}{g*_\Lb l}{glg^{-1}}$$
$$\fonctio{\Gb \times_\Pb \OCB.\Vb}{\Gb}{g*_\Pb x}{gxg^{-1}}\leqno{\mathit{and}}$$
are \'etale. 
\end{lem}

\bi

\proof By Remark \ref{identification}, 
it is sufficient to prove that the morphism 
$$\fonction{f}{\Gb \times_\Pb \OCB.\Vb}{\Gb}{g *_\Pb x}{gxg^{-1}}$$
is \'etale. 

Since $\Gb \times_\Pb \OCB.\Vb$ and $\Gb$ are smooth varieties, 
$f$ is \'etale if and only if the 
differential of $f$ at any point of $\Gb \times_\Pb \OCB.\Vb$ 
is an isomorphism \cite[Proposition III.10.4]{hartshorne}. 
By $\Gb$-equivariance of the morphism $f$ ($\Gb$ acts on 
$\Gb \times_\Pb \OCB.\Vb$ by left translation on the first factor, 
and acts on $\Gb$ by conjugation), it is sufficient to 
prove that $(df)_{1*_\Pb x}$ is an isomorphism for every 
$x \in \OCB.\Vb$. 

For this, let $\Pb^-$ denote the \para of $\Gb$ opposed to 
$\Pb$ (with respect to $\Lb$), and let $\Vb^-$ denote its unipotent 
radical. Then $\Vb^- \times\OCB.\Vb$ is an open neighborhood of 
$1 *_\Pb x$ in $\Gb *_\Pb \OCB.\Vb$. Therefore, it is sufficient 
to prove that the differential of the map 
$$\fonction{f^-}{\Vb^- \times \OCB.\Vb}{\Gb}{(g,l,x)}{glxg^{-1}}$$
at $(1,x)$ is an isomorphism for every $x \in \OCB.\Vb$.

Let $\gGB$, $\vGB^-$, $\lGB$ and $\pGB$ denote the Lie 
algebra of $\Gb$, $\Vb^-$, $\Lb$ and $\Pb$ respectively. 
Since $\OCB$ is open in $\Lb$, we may identify the tangent space 
to $\OCB.\Vb$ at $x$ with $\pGB$ (using left translation by $x$). 
Similarly, we identify the tangent space to $\Gb$ 
at $x$ with $\gGB$ using left translation. Using these identifications, 
the differential of $f^-$ at $(1,x)$ may be identified with 
the map
$$\fonction{\d}{\vGB^- \oplus \pGB}{\gGB}{A \oplus B }{
(\ad x)^{-1}(A) - A + B .}$$
For dimension reasons, we only need to prove that $\d$ is injective. 

For this, let $\l \in Y(\Tb)$ be such that $\Lb=C_\Gb(\Im \l)$ and 
$$\Pb=\{g \in \Gb~|~\lim_{t \to 0} \l(t)g\l(t)^{-1}~{\mathrm{exists}}\}.$$
For the definition of $\lim_{t \to 0} \l(t)$, see \cite[Page 184]{DLM1}. 
We then define, for each $i \in \ZM$, 
$$\gGB(i)=\{X \in \gGB~|~(\ad \l(t))(X)=t^i X\}.$$
Then
$$\gGB=\mathop{\oplus}_{i \in \ZM}~ \gGB(i),$$
$$\pGB=\mathop{\oplus}_{i \ge 0}~ \gGB(i),$$
$$\vGB^-=\mathop{\oplus}_{i < 0}~ \gGB(i)\leqno{\mathrm{and}}$$
(see \cite[5.14]{DLM1}). For each $X \in \gGB$, we denote by $X_i$ its projection 
on $\gGB(i)$.

Now, let $l=\pi_\Lb(x)$. Then it is 
clear that we have, for any $i_0 \in \ZM$ and any $X \in \gGB(i_0)$, 
$$(\ad x)^{-1}(X) \in (\ad l)^{-1}(X) + 
(\mathop{\oplus}_{i > i_0} \gGB(i)).\leqno{(1)}$$

Now, let $A \oplus B \in \Ker \d$, and assume that $A \not= 0$. 
Then there exists $i_0 < 0$ minimal among all $i < 0$ such that 
$A_i \not= 0$. Then, by equality (1), the projection of $\d(A \oplus B)$ on 
$\gGB(i_0)$ is equal to $(\ad l)^{-1}(A_{i_0})-A_{i_0}$. But, 
$\d(A \oplus B)=0$, so $(\ad l)^{-1}(A_{i_0})-A_{i_0}=0$. 
Therefore, $C_\gGB(l) \notincl \lGB$. So 
$C_\gGB(s) \notincl \lGB$, where $s$ denotes the semisimple part 
of $l$. 
However, $l$ lies in $\OCB$, so its semisimple part $s$ lies also 
in $\OCB$ by Lemma \ref{equiv class}. But, by \cite[Proposition 9.1 (1)]{borel}, 
$C_\gGB(s)$ is the Lie algebra of $C_\Gb^\ci(s)$ which 
is contained in $\Lb$ by Lemma \ref{equiv class}. We get a contradiction.

So, this discussion shows that $A=0$. But $0=\d(A,B)=(\ad x)^{-1}(A) - A + B$, 
so $B=0$. This completes the proof of Lemma \ref{lissable}.\fin

\bi

\sub{Isolated class} An element $g \in \Gb$ is said {\it ($\Gb$-)isolated} 
if the centralizer of its semisimple part is not contained in 
a \levi of a proper \para of $\Gb$.

Let $\Lb_\isol$ denote the subset of $\Lb$ consisting of 
$\Lb$-isolated elements, and let $\Tb_\isol=\Tb \cap \Lb_\isol$. 
Then $\Tb_\isol$ is a closed subset of $\Tb$. Therefore 
$\na_\Lb(\Tb_\isol)$ is a closed subset of $\Tb/W_\Lb$. 
Moreover, we have $\Lb_\isol=\na_\Lb^{-1}(\na(\Tb_\isol))$, 
so $\Lb_\isol$ is a closed subset of $\Lb$. 

As a consequence of Lemma \ref{proj mor}, the image of the morphism 
$\Gb \times_\Pb \Lb_\isol.\Vb \to \Gb$, $g *_\Pb x \mapsto gxg^{-1}$ is 
a closed subvariety of $\Gb$~: we denote it by $\Xb_{\Gb,\Lb}$. 
On the other hand, by Lemma \ref{lissable}, the image of the morphism 
$\Gb \times_\Lb \OCB \to \Gb$, $g *_\Lb x \mapsto gxg^{-1}$ 
is an open subset of $\Gb$, which will be denoted by $\OCB_{\Gb,\Lb}$. 
Finally, we denote by $\Yb_{\Gb,\Lb}$ the intersection of $\Xb_{\Gb,\Lb}$ 
and $\OCB_{\Gb,\Lb}$.

Note that $W_\Gb(\Lb)=N_\Gb(\Lb)/\Lb$ acts (on the right) on 
the variety $\Gb \times_\Lb \Lb_{\isol,\reg}$ in the following way. 
If $w \in W_\Gb(\Lb)$ and if $g *_\Lb l \in \Gb \times_\Lb \Lb_{\isol,\reg}$, 
then 
$$(g *_\Lb l).w=g\wdo *_\Lb \wdo^{-1} l \wdo,$$
where $\wdo \in N_\Gb(\Lb)$ is any representative of $w$.

\bi

\begin{prop}\label{etalisable}
The map $\Gb \times_\Lb \Lb_{\isol,\reg} \to \Yb_{\Gb,\Lb}$, 
$g *_\Lb l \mapsto glg^{-1}$ 
is a Galois \'etale covering with group $W_\Gb(\Lb)$.
\end{prop}

\bi

\proof Set $\pi : \Gb \times_\Lb \Lb_{\isol,\reg} \to \Yb_{\Gb,\Lb}$, 
$g *_\Lb l \mapsto glg^{-1}$, and let 
$\g : \Gb \times_\Lb \OCB \to \Gb$, $g *_\Lb l \mapsto glg^{-1}$. 
By Lemma \ref{lissable}, $\g$ is an \'etale morphism. If we prove 
that the square
$$\diagram 
\Gb \times_\Lb \Lb_{\isol,\reg} \rrto \ddto_{\DS{\pi}} && 
\Gb \times_\Lb \OCB \ddto^{\DS{\g}} \\
&&\\
\Yb_{\Gb,\Lb} \rrto && \Gb 
\enddiagram\leqno{(\#)}$$
is cartesian, then, by base change, we get that $\pi$ is an 
\'etale morphism. 

Since $\g$ is smooth, the fibred product (scheme) 
of $\Gb \times_\Lb \OCB$ and $\Yb_{\Gb,\Lb}$ over $\Gb$ 
is reduced (because $\Yb_{\Gb,\Lb}$ is), so we only need 
to prove that $\Gb \times_\Lb \Lb_{\isol,\reg} =\g^{-1}(\Yb_{\Gb,\Lb})$. 

So, let $g *_\Lb l \in \Gb \times_\Lb \OCB$ be such that 
$glg^{-1} \in \Yb_{\Gb,\Lb}$. Then there exist $h \in \Gb$, 
$m \in \Lb_\isol$ and $v \in \Vb$ such that $hmvh^{-1}=glg^{-1}$. 
Let $s$ (\resp $t$, \resp $t'$) denote the semisimple part of $l$, $m$ 
and $mv$. By Lemma \ref{V conjugue}, there exists $x \in \Vb$ 
such that $t'=\lexp{x}{t}$. Therefore $\lexp{hx}{t}=\lexp{g}{s}$. 

Since $t$ is $\Lb$-isolated, we have $\Zb(C_\Lb^\ci(t))^\ci=\Zb(\Lb)^\ci$. 
Therefore, $\Zb(C_\Gb^\ci(t))^\ci \incl \Zb(\Lb)^\ci$. 
On the other hand, since $s \in \OCB$, we have $C_\Gb^\ci(s) \incl \Lb$, 
so $\Zb(\Lb)^\ci \incl \Zb(C_\Gb^\ci(s))$. This proves that 
$\lexp{g}{\Zb(\Lb)^\ci} \incl \lexp{hx}{\Zb(\Lb)^\ci}$. For dimension 
reasons, we have $\lexp{g}{\Zb(\Lb)^\ci}=\lexp{hx}{\Zb(\Lb)^\ci}$,  
so $\Zb(\Lb)^\ci=\Zb(C_\Gb^\ci(s))$. Hence, $l$ is isolated. 
So, we have proved that $\pi$ is \'etale. 

Now $W_\Gb(\Lb)$ acts freely on $\Gb \times_\Lb \Lb_{\isol,\reg}$. 
So the quotient morphism 
$\Gb \times_\Lb \Lb_{\isol,\reg} \to \Gb \times_{N_\Gb(\Lb)} \Lb_{\isol,\reg}$ 
is a Galois \'etale covering with group $W_\Gb(\Lb)$. 
Moreover, $\pi$ clearly factorizes through this quotient 
morphism. We get an \'etale morphism 
$\pi_0 : \Gb \times_{N_\Gb(\Lb)} \Lb_{\isol,\reg} \to \Yb_{\Gb,\Lb}$. 
Proving Proposition \ref{etalisable} is now equivalent to prove 
that $\pi_0$ is an isomorphism of varieties. Since it is \'etale, 
we only need to prove that it is bijective. 

First $\pi_0$ is clearly surjective. We just need to prove that 
it is injective. Let $(g,l)$ and $(g',l')$ in $\Gb \times \Lb_{\isol,\reg}$ 
be such that $glg^{-1}=g'l'g^{\prime -1}$. Let $s$ and $s'$ be the semisimple 
parts of $l$ and $l'$ respectively. 
Then $\Zb(C_\Lb^\ci(s))^\ci=\Zb(\Lb)^\ci$ since $l$ is $\Lb$-isolated. 
Since $l \in \OCB$, we also have $C_\Gb^\ci(s)=C_\Lb^\ci(s)$, so 
$\Zb(C_\Gb^\ci(s))^\ci=\Zb(\Lb)^\ci$. Similarly, 
$\Zb(C_\Gb^\ci(s'))^\ci=\Zb(\Lb)^\ci$. So 
$\lexp{g^{-1}g'}{\Zb(\Lb)^\ci}=\Zb(\Lb)^\ci$. Therefore, $g^{-1}g' \in N_\Gb(\Lb)$. 
This completes the proof.\fin

\bi

\sec{Action of the relative Weyl group\label{sub relatif}}~ 

\med

\sub{The set-up} {\it From now on, and until the end of this paper, we 
denote by $\Sigb$ the inverse image of an $\Lb/\Zb(\Lb)^\ci$-isolated 
class of $\Lb/\Zb(\Lb)^\ci$ under the canonical projection 
$\Lb \to \Lb/\Zb(\Lb)^\ci$.} We also fix an element $v \in \Sigb$.


Following \cite[\SEC\SEC 3 and 4]{luicc}, we consider the varieties
$$\Ybh=\Gb \times \Sigb_\reg,$$
$$\Ybt=\Gb \times_\Lb \Sigb_\reg,$$
$$\Xbh=\Gb \times \overline{\Sigb}\Vb$$
$$\Xbt=\Gb \times_\Pb \overline{\Sigb}\Vb.\leqno{\mathrm{and}}$$
In these definitions, the group $\Lb$ (\resp $\Pb$) acts on $\Gb$ by right 
translations, and acts on $\Sigb_\reg$ (\resp $\overline{\Sigb} \Vb$) 
by conjugation. We also set 
$$\Yb=\bigcup_{g \in \Gb} g \Sigb_\reg g^{-1}$$
$$\Xb=\bigcup_{g \in \Gb} g \overline{\Sigb} \Vb g^{-1}.\leqno{\mathrm{and}}$$
By Lemma \ref{V conjugue}, we have 
\equat\label{x}
\Xb \incl \na^{-1}\bigl(\na(\Zb(\Lb)^\ci)\bigr)
\endequat
and
\equat\label{y}
\Yb \incl \na^{-1}\bigl(\na(\Zb(\Lb)^\ci_\reg)\bigr).
\endequat
Moreover, $\Xb$ is the image of $\Gb \times_\Pb \overline{\Sigb}.\Vb$ under 
the morphism $\Gb \times_\Pb \Pb \to \Gb$, $g *_\Pb x \mapsto gxg^{-1}$ 
studied in Subsection \ref{sub fam}. So, by Lemma \ref{proj mor}, 
$\Xb$ is a closed irreducible subvariety of $\Gb$. We set 
$$\fonction{\piba}{\Gb \times_\Pb \overline{\Sigb}.\Vb}{\Xb}{g*_\Pb x}{gxg^{-1}.}$$
It is a projective morphism of varieties.

On the other hand, $\Yb^+=\bigcup_{g \in \Gb} g (\overline{\Sigb})_\reg g^{-1}$ 
is in fact the intersection of $\Xb$ with $\na^{-1}\bigl(\na(\Zb(\Lb)^\ci_\reg)\bigr)$, 
so it is an open subset of $\Xb$ (by Lemma \ref{nabla -2}). But, by 
Proposition \ref{etalisable}, 
the map $\Gb \times_\Lb (\overline{\Sigb})_\reg \to \Yb^+$, 
$g *_\Lb l \mapsto glg^{-1}$ is a Galois 
\'etale covering with group $W_\Gb(\Lb,\Sigb)=N_\Gb(\Lb,\Sigb)/\Lb$. 
Since $\Gb \times_\Lb \Sigb_\reg$ is open in $\Gb \times_\Lb (\overline{\Sigb})_\reg$, 
its image $\Yb$ under this \'etale morphism is open in $\Yb^+$. 
This proves that $\Yb$ is open in $\Xb$. Moreover, 
since the map 
$$\fonction{\pi}{\Gb \times_\Lb \Sigb_\reg}{\Yb}{g *_\Lb l}{glg^{-1}}$$
is a Galois \'etale covering with group $W_\Gb(\Lb,\Sigb)$, 
we get that $\Yb$ is smooth (indeed, $\Gb \times_\Lb \Sigb_\reg$ is smooth). 

Recall that $\Gb \times_\Lb \Sigb_\reg \to \Gb \times_\Pb \Sigb_\reg.\Vb$, 
$g *_\Lb l \mapsto g*_\Pb l$ is an isomorphism (see Corollary \ref{U conjugue bis}). 
Moreover, it is clear that $\piba^{-1}(\Yb)=\Gb \times_\Pb \Sigb_\reg.\Vb$. 
We summarize all these facts in the following proposition.

\bi

\begin{prop}[{\bf Lusztig \cite[\P 3.1, 3.2 and Lemma 4.3]{luicc}}]\label{resume truc}
With the above notation, we have~:

\tete{1} $\Xb$ is a closed irreducible subvariety of $\Gb$ and 
$\Yb$ is open in $\Xb$.

\tete{2} The natural map $\Ybt \to \Xbt$, $g *_\Lb x \mapsto g *_\Pb x$ is an 
open immersion and the square
$$\diagram 
\Ybt \rrto^{\DS{\pi}} \ddto && \Yb \ddto\\
&&\\
\Xbt 
\rrto^{\DS{\piba}}  && \Xb
\enddiagram$$
is cartesian.

\tete{3} $\piba$ (hence $\pi$) is a projective morphism.

\tete{4} $\pi$ is an \'etale Galois covering with group $W_\Gb(\Lb,\Sigb)$. 
In particular, $\Yb$ is smooth.
\end{prop}

\bi

Note that $\Gb$ acts on $\Ybh$, $\Ybt$, $\Xbh$ and $\Xbt$ 
by left translation on the first factor, and that it acts on $\Yb$ and $\Xb$ 
by conjugation. Also, the group $\Zb(\Gb) \cap \Zb(\Lb)^\ci$ acts 
on the varieties $\Ybh$, $\Ybt$, $\Xbh$ and $\Xbt$ 
by left translation on the second factor, and it acts on $\Yb$ and $\Xb$ 
by left translation. These actions of $\Gb$ and $\Zb(\Gb) \cap \Zb(\Lb)^\ci$ 
commute.

We have the following commutative diagram
\equat\label{diagramme central}
\diagram
\Sigb_\reg \ddto && \Ybh \llto_{\DS{\a}} 
\rrto^{\DS{\b}} \ddto && \Ybt \rrto^{\DS{\pi}} \ddto && \Yb \ddto\\
&&&&&&\\
\overline{\Sigb} && \Xbh \llto_{\DS{\alpba}} \rrto^{\DS{\betba}} && \Xbt 
\rrto^{\DS{\piba}}  && \Xb,
\enddiagram
\endequat
where $\a$ and $\alpba$ are the canonical projections, and $\b$ and $\betba$ 
are the canonical quotient morphisms
Moreover, the vertical maps are the natural ones. 
Note that the morphisms $\b$, $\betba$, $\pi$ and $\piba$ are $\Gb \times 
(\Zb(\Gb) \cap \Zb(\Lb)^\ci)$-equivariant, 
and that the vertical maps $\Ybh \to \Xbh$, $\Ybt \to \Xbt$ and $\Yb \to \Xb$ 
are also $\Gb \times (\Zb(\Gb) \cap \Zb(\Lb)^\ci)$-equivariant.

\bi

\sub{Extension of the action of $W_\Gb(\Lb,\Sigb)$} 
We must notice that the action of the group $W_\Gb(\Lb,\Sigb)$ is defined only 
on $\Ybt$. 
However, we will see in this subsection that it is possible to 
extend it to an open subset $\Xbt_\mini$ of $\Xbt$ 
which, in general, contains strictly $\Ybt$.

We first need some preliminaries to construct this extension. 
If $l \in \overline{\Sigb}$, then 
$$\dim C_\Lb(l) \ge \dim C_\Lb(v)$$
and equality holds if 
and only if $l \in \Sigb$. Consequently, if $g \in \Xb$, then 
$$\dim C_\Gb(g) \ge \dim C_\Lb(v)$$ 
(\cf Lemma \ref{dim p} (2)). Moreover, if 
the equality holds, then, by Lemma \ref{dim p} and the previous 
remark, $\piba^{-1}(g)$ is contained in 
$$\Xbt_0=\Gb \times_\Pb \Sigb.\Vb,$$
which is a smooth open subset of $\Xbt$.

Let us define now
$$\Xb_\mini=\{g \in \Xb~|~ \dim C_\Gb(g) = \dim C_\Lb(v)\}$$
$$\Xbt_\mini=\piba^{-1}(\Xb_\mini).\leqno{\mathrm{and}}$$    
By \cite[Proposition 1.4]{humphreys}, $\Xb_\mini$ is an open subset of $\Xb$, and, 
by the previous discussion, 
$\Xbt_\mini \incl \Xbt_0$. Also, $\Yb \incl \Xb_\mini$, so $\Ybt \incl \Xbt_\mini$.
Now, let $\pi_\mini : \Xbt_\mini \to \Xb_\mini$ 
denote the restriction of $\piba$~: it is a projective morphism. 

Moreover, if $g \in \Xb_\mini$, then, by Lemma \ref{dim p} (1), 
$\pi_\mini^{-1}(g)$ is a finite set. 
The \mor $\pi_{\min}$ being projective and quasi-finite, it is finite 
\cite[Exercise III.11.2]{hartshorne}. We gather these facts in the 
next proposition.

\bi

\begin{prop}\label{tas}
With the above notation, we have~:

\tete{1} $\Xb_\mini$ is a $\Gb \times (\Zb(\Gb) \cap \Zb(\Lb)^\ci)$-stable 
open subset of $\Xb$ containing $\Yb$.

\tete{2} $\Xbt_\mini$ is a $\Gb \times (\Zb(\Gb) \cap \Zb(\Lb)^\ci)$-stable 
smooth open subset of $\Xbt$ containing $\Ybt$.

\tete{3} The \mor $\pi_\mini : \Xbt_\mini \to \Xb_\mini$ is finite.
\end{prop}

\bi

The Proposition \ref{tas} has the following immediate consequence~:

\bi

\begin{theo}\label{extension action}
$(a)$ The variety $\Xbt_\mini$ is the normalization of $\Xb_\mini$ in the variety 
$\Ybt$. Therefore, there exists a unique action of the finite group $W_\Gb(\Lb,\Sigb)$ 
on the variety $\Xbt_\mini$ extending its action on $\Ybt$. 

\tete{b} This action is $\Gb \times (\Zb(\Gb) \cap \Zb(\Lb)^\ci)$-equivariant, and 
the \mor $\pi_\mini$ factorizes through the quotient $\Xbt_\mini/W_\Gb(\Lb,\Sigb)$.

\tete{c} If the variety $\Xb_\mini$ is normal, then $\pi_\mini$ induces an \iso 
of varieties $\Xbt_\mini/W_\Gb(\Lb,\Sigb) \simeq \Xb_\mini$.
\end{theo}

\bi

\noindent{\sc Notation -} (1) Let $(\Sigb.\Vb)_\mini$ denote the open subset 
$\Sigb.\Vb \cap \Xb_\mini$ of $\Sigb.\Vb$. Then $\Xbt_\mini=\Gb \times_\Pb 
(\Sigb.\Vb)_\mini$. 

\med

\tete{2} If there is some ambiguity, we will denote by $?_\Lb^\Gb$ the 
object $?$ defined above (for instance, $\Ybh_\Lb^\Gb$, $\Xbt_{\min,\Lb}^\Gb$, 
$\Xb_\Lb^\Gb$, $\pi_{\min,\Lb}^\Gb$...).


\bi

\sub{Unipotent classes} {\it From now on, and until the end of this paper, 
$\Sigb$ is the inverse image of a unipotent class of $\Lb/\Zb(\Lb)^\ci$.} 
Note that a unipotent class is isolated. Let $\Cb$ denote 
the unique unipotent class contained in $\Sigb$. From now on, the 
element $v$ introduced in the previous section is chosen 
in $\Cb$. Note that $\Sigb=\Zb(\Lb)^\ci.\Cb \simeq \Zb(\Lb)^\ci \times \Cb$ 
and that $\Sigb_\reg=\Zb(\Lb)^\ci_\reg.\Cb$.

\bi

\noindent{\sc Notation - } 
(1) We denote by $\Cb^\Gb$ the induced unipotent class of $\Cb$ 
from $\Lb$ to $\Gb$, that is the unique unipotent class $\Cb_0$ of $\Gb$ 
such that $\Cb_0 \cap \Cb.\Vb$ is dense in $\Cb.\Vb$.

\tete{2} If $z \in \Zb(\Lb)^\ci$, the group $C_\Pb(z)=\Lb.C_\Vb(z)$ is connected 
(\cf Lemma \ref{spa connexe}). So it is a \para of $C_\Gb^\ci(z)$ 
by \cite[Proposition 1.11 (ii)]{dmgnc}, with unipotent radical $C_\Vb(z)$ and 
Levi factor 
$\Lb$. We denote by $u_z$ an element of 
$\Cb^{C_\Gb^\ci(z)} \cap vC_\Vb(z)$. We set $\uti_z=1*_\Pb zu_z 
\in \Xbt_\mini$. 

\med

\tete{3} For simplicity, the unipotent element $u_1$ will be denoted by $u$, 
and $\uti$ stands for $\uti_1$. 

\bi

\remark{description xmin} Let us investigate here what are the elements 
of $\Xb_{\min}$. Since $\Xbt_\mini \incl \Xbt_0$, we need only to determine 
which elements of $\Sigb.\Vb$ belong to $\Xb_{\min}$. 
Let $g \in \Sigb.\Vb$. Let $z$ (\resp $u'$) be the semisimple 
(\resp unipotent) part of $g$. By Lemma \ref{V conjugue}, we may assume 
that $z$ belongs to $\Zb(\Lb)^\ci$. Now, let $\Gb'=C_\Gb^\ci(z)$, 
$\Pb'=C_{\Pb}^\ci(z)$, and $\Vb'=C_\Vb(z)$. Then $\Gb'$ is a reductive subgroup 
of $\Gb$ containing $\Lb$, $\Pb'$ is a \para of $\Gb'$, and $\Vb'$ is its 
unipotent radical. Then, by 
\cite[Corollary 11.12]{borel}, we have $u' \in \Gb'$. On the other hand, 
$C_\Gb^\ci(g) = C_{\Gb'}^\ci(u')$. Now, by Lemma \ref{dim p} (2) and by 
\cite[Proposition II.3.2 (b) and (e)]{spaltenstein}, $g \in \Xb_{\min}$ if and only 
if $u' \in \Cb^{\Gb'}$.

Hence, we have proved that
\equat\label{une seule}
\na^{-1}(\na(z)) \cap \Xb_\mini =(zu_z)_\Gb
\endequat
for every $z \in \Zb(\Lb)^\ci$, and that
$$\Xb_\mini=\bigcup_{z \in \Zb(\Lb)^\ci} (zu_z)_\Gb.~\SS{\square}$$

\bi

If $z \in \Zb(\Lb)^\ci$, 
we denote by $H_\Gb(\Lb,\Sigb,z)$ the stabilizer 
of $\uti_z$ in $W_\Gb(\Lb,\Sigb)$. 
We first investigate what is the group $H_\Gb(\Lb,\Sigb,1)$. 
Recall that the group $C_\Gb^\ci(u)$ is 
contained in $\Pb$ \cite[Proposition II.3.2 (e)]{spaltenstein}, so that 
$C_\Gb(u)/C_\Pb(u)$ is a finite set.

\bi

\begin{lem}\label{agu general}
We have $\pi_\mini^{-1}(u)=\{g *_\Pb u~|~g \in C_\Gb(u)\}$. In particular, 
$$|\pi_\mini^{-1}(u)|=|C_\Gb(u)/C_\Pb(u)|=|A_\Gb(u)/A_\Pb(u)|.$$
\end{lem}

\bi

\proof It is clear that $\{g *_\Pb u~|~g \in C_\Gb(u)\}$ is contained in 
$\pi_\mini^{-1}(u)$. Conversely, let $g *_\Pb x \in \pi_\mini^{-1}(u)$. 
By replacing $(g,x)$ by $(gl^{-1},lx)$ for a suitable choice of 
$l \in \Lb$, we may assume that $\pi_\Lb(x)=v$. Since $gxg^{-1} =u$, 
this means that $x \in v\Vb \cap \Cb^\Gb$. By 
\cite[Proposition II.3.2 (d)]{spaltenstein}, there exists $y \in \Pb$ such 
that $yxy^{-1}=u$. Therefore, $gy^{-1} \in C_\Gb(u)$ and 
$g*_\Pb x= gy^{-1} *_\Pb u$.\fin

\bi

\begin{coro}\label{agu 1}
If $C_\Gb(u) \incl \Pb$, then $\pi_\mini^{-1}(u)=\{\uti\}$. In particular, 
$W_\Gb(\Lb,\Sigb)$ stabilizes $\uti$, that is $H_\Gb(\Lb,\Sigb,1)=W_\Gb(\Lb,\Sigb)$. 
\end{coro}

\bi

Now, let us consider the general case. The second projection $\Xbh
\simeq \Gb \times \Zb(\Lb)^\ci \times \overline{\Cb} \times \Vb \to \Zb(\Lb)^\ci$ 
factorizes through the quotient morphism $\Xbh \to \Xbt$. We 
denote by $\zb : \Xbt \to \Zb(\Lb)^\ci$ the morphism obtained after 
factorization. The group $W_\Gb(\Lb)$ acts on $\Zb(\Lb)^\ci$ by conjugation, 
and it is easy to check that the restriction $\zb_\reg : \Ybt \to \Zb(\Lb)^\ci_\reg$ 
of $\zb$ to $\Ybt$ is $W_\Gb(\Lb,\Sigb)$-equivariant. Hence, 
the morphism $\zb_\mini : \Xbt_\mini \to \Zb(\Lb)^\ci$ obtained by 
restriction from $\zb$ is $W_\Gb(\Lb,\Sigb)$-equivariant.

As a consequence, we get 
\equat\label{ZZZ}
\Stab_{W_\Gb(\Lb,\Sigb)}(\gti) \incl \Stab_{W_\Gb(\Lb,\Sigb)}(\zb(\gti))
\endequat
for every $\gti \in \Xbt_\mini$. Also, note that $\zb(\gti)$ is conjugate 
in $\Gb$ to the semisimple part of $g=\pi_\mini(\gti)$ (\cf Proposition 
\ref{V conjugue}).
In fact, one can easily get a better result~:

\bi

\begin{prop}\label{stabilisateur z}
Let $z \in \Zb(\Lb)^\ci$. Then~:

\tete{1} $H_\Gb(\Lb,\Sigb,z)=H_{C_\Gb^\ci(z)}(\Lb,\Sigb,1)$.


\tete{2} If $C_{C_\Gb^\ci(z)}(u_z) \incl \Pb$, then 
$H_\Gb(\Lb,\Sigb,z) = W_{C_\Gb^\ci(z)}(\Lb,\Sigb)$.
\end{prop}

\bi

\proof For the proof of (1), the reader may refer to the proof of 
Proposition \ref{z} of this paper~: indeed, the situations are quite 
analogous, and the arguments involved are exactly similar. However, 
since the situation in Proposition \ref{z} is a little more complicated,
we have decided to give a detailed proof only in this case. 
(2) follows from (1) and from Corollary \ref{agu 1}.\fin

\bi

\sub{An example} 
Assume in this subsection, and only in this subsection, that $\Lb=\Tb$. 
Then $\Cb=1$, $\Sigb=\Tb$, $\Xbt = \Gb \times_\Bb \Bb$, 
$\Xb=\Gb$ and $\piba : \Xbt \to \Gb$ is the well-known Grothendieck map. Also, 
$W_\Gb(\Lb,\Sigb)=W$ in this case. 
Moreover, $\Xb_\mini$ is the open subset of $\Gb$ consisting of regular elements. 
As an open subset of $\Gb$, it is smooth. So the action 
of $W$ on $\Ybt$ extends to $\Xbt_\mini$ and $\Xbt_\mini/W = \Xb_\mini$.

Now, let $\gti \in \Xbt_\mini$, $g=\pi_\mini(\gti)$, and $t=\zb(\gti) \in \Tb$. 
We denote by $W^\ci(t)$ the Weyl group of $C_\Gb^\ci(t)$ relative to $\Tb$. The fiber 
$\pi_\mini^{-1}(g)$ may be identified with the set of Borel subgroups 
of $\Gb$ containing $g$. Since $\Xbt_\mini/W = \Xb_\mini$, $W$ acts transitively 
on $\pi_\mini^{-1}(g)$. But, $u_t$ is a regular unipotent element of $C_\Gb^\ci(t)$. 
Therefore, $C_{C_\Gb^\ci(t)}(u_t) \incl \Bb$. So, by Proposition \ref{stabilisateur z}, 
we have
$$\Stab_W(\gti)=W^\ci(t).$$
As a consequence, we get the well-known result 
\equat\label{nombre borel}
|\{x\Bb \in \Gb/\Bb~|~g \in \lexp{x}{\Bb}\}|=|W|/|W^\ci(t)|.
\endequat

\bi

\rem It is not the easiest way to prove \ref{nombre borel}~!

\bi

{\small
\example{sl2 2} Assume in this example that $\Gb=\Gb\Lb_2(\FM)$, that 
$$\Lb=\Tb=\{\diag(a,b)~|~a,b \in \FM^\times\},$$ 
and that 
$\Sigb=\Tb$. Let us denote by $\Pb^1$ the projective line. Then
$$\Xbt\simeq\{(\matrice{a & b \\ c & d},[x,y]) \in \Gb \times \Pb^1~|~ 
[ax+by,cx+dy]=[x,y]\},$$
$$\Xb=\Gb,$$ 
and $\pi : \Xbt \to \Xb$ is identified with the first projection. 
Moreover, $\Xb_\mini$ is the open subset of $\Gb$ consisting of 
non-central elements. We shall give a precise formula for describing 
the action of $W$ on $\Xbt_\mini$ in this little example.

Let $w$ denote the unique non-trivial element of $W$. 
It has order $2$. We define the right action of $w$ on 
$(g,[x,y]) \in \Xbt_\mini$ by
$$(g,[x,y]).w=\left\{\begin{array}{ll}
(g,[bx,(d-a)x-by]) & \quad{\mathrm{if}}~(bx,(d-a)x-by) 
\not=(0,0) \\
&\\
(g,[(a-d)y-cx,cy]) & \quad{\mathrm{if}}~((a-d)y-cx,cy)\not=(0,0),
\end{array}\right.$$
where $g=\matrice{a & b \\ c & d}$. 
One can check that, if $\Xbt_1$ (\resp $\Xbt_2$) is the open subset 
of $\Xbt$ defined by the first condition (\resp the second condition), 
then $\Xbt_1 \cup \Xbt_2 = \Xbt_\mini$, and that the formulas given above 
coincide on $\Xbt_1 \cap \Xbt_2$. So we have defined a \mor of varieties. 
The fact that it is an \auto of order $2$ is obvious, and the reader can 
check that it extends the action of $W$ on $\Ybt$. 

One can also check, as it is expected by \ref{nombre borel}, 
that $W$ acts trivially on the elements 
$(g,[x,y]) \in \Xbt_\mini$ such that $g$ is not semisimple.\finl}

\bi

\sec{A morphism $W_\Gb(\Lb,\Sigb) \to A_\Lb(v)$\label{mor phi}}~ 

\med

The restriction of an $\Lb$-equivariant local system on $\Cb$ 
through the morphism $\Lb/C_\Lb^\ci(v) \to \Cb$, 
$lC_\Lb^\ci(v) \mapsto lvl^{-1}$ is constant. That is why this 
morphism is interesting when one is working with character 
sheaves (which are equivariant intersection cohomology complexes). 
This morphism 
can be followed all along the diagram \ref{diagramme central}, 
and it gives rise to new varieties on which the group 
$W_\Gb(\Lb,v)=N_\Gb(\Lb,v)/C_\Lb^\ci(v)$ acts 
(note that $W_\Gb(\Lb,v)/A_\Lb(v) \simeq W_\Gb(\Lb,\Sigb)$). 
Following the method of the previous section, these actions 
can be extended to some variety $\Xbt_\mini^\prime$ lying over $\Xbt_\mini$. 
We show in this section that some stabilizers under this action
can be described in terms of a morphism of groups 
$W_\Gb(\Lb,\Sigb) \to A_\Lb(v)$ 
(under some little hypothesis). In \SEC\ref{ele sec}, elementary 
properties of this morphism will be investigated.

\bi

\sub{Notation} Let $\Sigb'=\Lb/C_\Lb^\ci(v) \times \Zb(\Lb)^\ci$, and 
let $\Sigb_\reg^\pr = \Lb/C_\Lb^\ci(v) \times \Zb(\Lb)_\reg^\ci $. We denote by $f : 
\Sigb' \to \Sigb$, $(lC_\Lb^\ci(v),z) \mapsto lzvl^{-1}=zlvl^{-1}$. Then $f$ is a 
finite surjective $\Lb$-equivariant morphism (here, $\Lb$ acts 
on $\Sigb'$ by left translation on the first factor). We denote by $f_\reg : 
\Sigb_\reg^\pr \to \Sigb_\reg$ the restriction of $f$. 

Now, let 
$$\Ybh'= \Gb \times \Sigb_\reg^\pr,$$
$$\Ybt'=\Gb \times_\Lb \Sigb_\reg^\pr = \Gb/C_\Lb^\ci(v) \times \Zb(\Lb)^\ci_\reg.
\leqno{\mathrm{and}}$$
We then get a commutative diagram
\equat\label{diagramme}
\diagram
\Sigb_\reg^\pr \ddto_{\DS{f_\reg}} && \Ybh' \llto_{\DS{\a'}} \rrto^{\DS{\b'}} 
\ddto_{\DS{\fha}} &&
 \Ybt' \ddto_{\DS{\fti}} \ddrrto^{\DS{\pi'}} & \\
&&&&&&\\
\Sigb_\reg && \Ybh \llto^{\DS{\a}} \rrto_{\DS{\b}} && \Ybt \rrto_{\DS{\pi}} && \Yb,
\enddiagram
\endequat
where the vertical maps are induced by $f_\reg$, $\a'$ is the projection 
on the first factor, $\b'$ is the quotient morphism, and $\pi'=\pi \ci \fti$. 
The group $\Gb$ acts on $\Ybh'$ and $\Ybt'$ by left translation on the 
first factor, and acts on $\Yb$ by conjugation. The group $\Zb(\Gb) \cap \Zb(\Lb)^\ci$ 
acts on $\Sigb_\reg^\pr$ by translation on the second factor~: it induces 
an action on $\Ybh^\pr$ and $\Ybt'$. The morphisms $\fha$, $\fti$, 
$\b'$ and $\pi'$ are $\Gb \times (\Zb(\Gb) \cap \Zb(\Lb)^\ci)$-equivariant. 
Moreover, all the squares of diagram \ref{diagramme} are cartesian.

Now, we define  
$$W_\Gb(\Lb,v) = N_\Gb(\Lb,v)/C_\Lb^\ci(v)$$ 
(note that $N_\Gb(\Lb,v)^\ci=C_\Lb^\ci(v)$).
The group $C_\Lb(v)$ is a normal subgroup of $N_\Gb(\Lb,v)$ so $A_\Lb(v)$ is a normal 
subgroup of $W_\Gb(\Lb,v)$. Note that
$$W_\Gb(\Lb,v)/A_\Lb(v) \simeq W_\Gb(\Lb,\Sigb).$$
The group $N_\Gb(\Lb,v)$ acts freely on the right 
on the variety $\Ybh'$ in the following way~: if $w \in N_\Gb(\Lb,v)$ 
and if $(g,lC_\Lb^\ci(v),z) \in \Ybh'$, then
$$(g,lC_\Lb^\ci(v),z).w=(gw,w^{-1}lwC_\Lb^\ci(v),w^{-1}zw).$$
This induces a free right $\Gb \times (\Zb(\Gb) \cap \Zb(\Lb)^\ci)$-equivariant 
action of $W_\Gb(\Lb,v)$ on $\Ybt'$. 
Moreover, the fibers of the \mor $\pi'$ are $W_\Gb(\Lb,v)$-orbits.

\bi

\remark{action A} If $a \in A_\Lb(v)$ and $g *_\Lb (lC_\Lb^\ci(v),z) \in \Ybt'$, then 
$$(g *_\Lb (lC_\Lb^\ci(v),z)).a=g*_\Lb (laC_\Lb^\ci(c),z).~\SS{\square}$$

\bi

\sub{Normalization\label{sub nor}} 
Let $\Xbt'$ be the normalization of the variety $\Xbt$ in $\Ybt'$. We denote by 
$\fba : \Xbt' \to \Xbt$ the corresponding \mor of varieties. Let 
$\Xbt_0^\pr$ (\resp $\Xbt_\mini^\pr$) denote the inverse image, in $\Xbt'$, 
of the variety $\Xbt_0$ (\resp $\Xbt_\mini$). We denote by $\fti_0 : 
\Xbt_0^\pr \to \Xbt_0$ (\resp $\fti_\mini : \Xbt_\mini^\pr \to \Xbt_\mini$) 
the restriction of $\fba$ to $\Xbt_0^\pr$ (\resp $\Xbt_\mini^\pr$). 
Then $\Xbt_0^\pr$ (\resp $\Xbt_\mini^\pr$) is the normalization of $\Xbt_0$ 
(\resp $\Xbt_\mini$) in $\Ybt^\pr$. 
We can summarize the notation in the following commutative diagram
$$\diagram
\Ybt' \rrto \ddto_{\DS{\fti}}\xto'[2,1]^{\DS{\pi'}}[4,2]&& 
\Xbt_\mini^\pr \ddto_{\DS{\fti_\mini}}\xto'[2,1]^{\DS{\pi_\mini^\pr}}[4,2]
\rrto && \Xbt_0^\pr \ddto^{\DS{\fti_0}} \rrto && \Xbt' \ddto^{\DS{\fba}} 
\xto[4,2]^{\DS{\piba}'}\\
&&\\
\Ybt \xto[0,2] \ddrrto_{\DS{\pi}} && \Xbt_\mini \ddrrto_{\DS{\pi_\mini}} \rrto 
&& \Xbt_0  \rrto && \Xbt \ddrrto^{\DS{\piba}}\\
&&&&\\
&& \Yb \rrto && \Xb_\mini \xto[0,4] &&&&\Xb. \\
\enddiagram$$
In this diagram, all the horizontal maps are open immersions, 
and all the squares are cartesian. 
Since $\Xbt_\mini$ is the normalization of $\Xb_\mini$ in $\Ybt$, we get~:

\begin{theo}\label{tango}
We have~:

\tete{1} The variety $\Xbt_\mini^\pr$ is the normalization of $\Xb_\mini$ in $\Ybt^\pr$. 
Therefore, the action of $W_\Gb(\Lb,v)$ on $\Ybt'$ extends uniquely to 
an action of $W_\Gb(\Lb,v)$ on $\Xbt_\mini^\pr$.

\tete{2} $\Xbt'$ inherits from $\Ybt'$ an action of $\Gb \times (\Zb(\Gb) 
\cap \Zb(\Lb)^\ci)$, and this action commutes with the one of $W_\Gb(\Lb,v)$ 
on $\Xbt_\mini$.
\end{theo}

\bi

\remark{xtilde 0} We do not know how to determine in general the variety 
$\Xbt'$. However, it is possible to give an explicit description of $\Xbt_0^\pr$. 
This can be done as follows. The \para $\Pb$ acts on $\Sigb' \times \Vb$ 
by the following action~: 
if $l,l_0 \in \Lb$, $x, x_0 \in \Vb$, and $z_0 \in \Zb(\Lb)^\ci$, then 
$$\lexp{lx}{(l_0 C_\Lb^\ci(v),z_0,x_0)}=(ll_0C_\Lb^\ci(v),z_0,
\lexp{ll_0z_0vl_0^{-1}}{x}x_0 x^{-1}).$$
The reader can check that this defines an action of $\Pb$. Moreover, the 
\mor $f \times \Id_\Vb : 
\Sigb' \times \Vb \to \Sigb\Vb$, $(l_0C_\Lb^\ci(v),z_0,x_0) \mapsto 
l_0z_0vl_0^{-1}x_0$ induced by $f$ is $\Pb$-equivariant. 

By Corollary \ref{U conjugue bis}, 
$\Gb \times_\Lb \Sigb_\reg^\pr \simeq \Gb \times_\Pb (\Sigb_\reg^\pr \times \Vb)$ 
is an open subset of $\Gb \times_\Pb (\Sigb' \times \Vb)$ isomorphic 
to $\Ybt'$. Moreover, the morphism 
$$\Gb \times_\Pb (\Sigb' \times \Vb) \longto \Xbt=\Gb \times_\Pb \Sigb\Vb$$
induced by $f$ is finite, as it can be checked by restriction to an open subset 
of the form $g\Vb^- \Pb \times_\Pb (\Sigb' \times \Vb)$. Here, $\Vb^-$ 
is the unipotent radical of the opposite \para $\Pb^-$ of $\Pb$ 
with respect to $\Lb$. Finally, by the same argument, 
$\Gb \times_\Pb (\Sigb' \times \Vb)$ is smooth. Hence
\equat\label{normalization x0}
\Xbt_0^\pr=\Gb \times_\Pb (\Sigb' \times \Vb).
\endequat
Since $\Xbt'$ is the normalization of $\Xbt$ in $\Ybt'$, it inherits an action 
of the group $A_\Lb(v)$. It is very easy to describe this action on $\Xbt_0^\pr$ 
by using \ref{normalization x0}.
It acts on $\Xbt_0^\pr$ by right translation on the 
factor $\Lb/C_\Lb^\ci(v)$ of $\Sigb'$. This is a free action and 
the fibers of $\fti_0$ are $A_\Lb(v)$-orbits.\finl

\bi

We will denote by $(\Sigb' \times \Vb)_\mini$ the inverse image, under 
$f \times \Id_\Vb$, of the open subset $(\Sigb\Vb)_\mini$ of $\Sigb\Vb$. 
Then $\Xbt_\mini^\pr=\Gb \times_\Pb (\Sigb' \times \Vb)_\mini$. The action 
of $W_\Gb(\Lb,v)$ is quite mysterious, but the action of its subgroup $A_\Lb(v)$ 
is understandable. It is obtained by restriction from its action on $\Xbt_0^\pr$ 
which is described at the end of Remark \ref{xtilde 0}.

\bi

\sub{Stabilizers\label{02}} If $z \in \Zb(\Lb)^\ci$, let $\uti_z^\pr=
1 *_\Pb (C_\Lb^\ci(v),z,v^{-1}u_z) \in \Xbt_\mini^\pr$. Recall 
that $u_z$ is an element of $v C_\Vb(z) \cap \Cb^{C_\Gb^\ci(z)}$. Note that 
$\fti_\mini^\pr(\uti^\pr_z)=\uti_z$, so that 
$\pi_\mini^\pr(\uti^\pr_z)=u_z$. For simplification, we denote by $\uti'$ the element 
$\uti_1^\pr$. The stabilizer of the element $\uti_z^\pr$ in $W_\Gb(\Lb,v)$ is denoted 
by $H_\Gb(\Lb,v,z)$. The aim of this subsection is to get some informations 
about these stabilizers.

The first result comes from the fact that $A_\Lb(v)$ acts freely on $\Xbt_\mini^\pr$~:
\equat\label{H cap A}
H_\Gb(\Lb,v,z) \cap A_\Lb(v) =\{1\}.
\endequat
The second one is analogous to Proposition \ref{stabilisateur z}~: it may be viewed 
as a kind of Jordan decomposition.

\bi

\begin{prop}\label{z}
If $z \in \Zb(\Lb)^\ci$, then $H_\Gb(\Lb,v,z)=H_{C_\Gb^\ci(z)}(\Lb,v,1)$.
\end{prop}

\bi

\proof Let $\zb'$ denote the composite morphism of varieties 
$\Xbt' \longmapright{\fba} \Xbt \longmapright{\zb} \Zb(\Lb)^\ci$, and let 
$\zb_\mini^\pr : \Xbt_\mini^\pr \to \Zb(\Lb)^\ci$ denote the restriction of 
$\zb'$. Then $\zb_\mini^\pr$ is a $W_\Gb(\Lb,v)$-equivariant morphism (as it can be 
checked by restriction to $\Ybt'$). So, the group $H_\Gb(\Lb,v,z)$ is contained 
in $\WC_z=W_{C_\Gb(z)}(\Lb,v)$. 

Let $\Ab_z=\{t \in \Zb(\Lb)^\ci~|~C_\Gb^\ci(t) \incl C_\Gb^\ci(z)\}$~: 
it is an open subset of $\Zb(\Lb)^\ci$ containing $z$ and $\Zb(\Lb)^\ci_\reg$. 
Now let $\Sigb_z=\Ab_z .\Cb$ and let $\Sigb_z^\pr=\Lb/C_\Lb^\ci(v) \times \Ab_z$. 
Then 
$$\Xbt_z^\pr=\Gb\times_\Pb (\Sigb_z^\pr \times \Vb)_{\min,\Lb}^\Gb$$ 
is an open subset of $\Xbt_\mini$ containing $\uti_z$, and it is stable under 
the action of $\WC_z$, since $\Ab_z$ is and since 
$\Xbt_z^\pr=\zb_\mini^{\pr -1}(\Ab_z)$. Now, let 
$$\Xbt'(z)=
C_\Gb(z) \times_{C_\Pb(z)} (\Sigb_z^\pr \times C_\Vb(z))_{\mini,\Lb}^{C_\Gb^\ci(z)}.$$
The natural morphism $\Xbt'(z) \to \Xbt_z^\pr$ is injective and $\WC_z$-equivariant. 
This proves that the stabilizer $H_\Gb(\Lb,v,z)$ is equal to the 
stabilizer of $1 *_{C_\Pb(z)} (C_\Lb^\ci(v),z,v^{-1}u_z) \in \Xbt'(z)$ 
in $\WC_z$. But this stabilizer must stabilize the connected component 
of $1 *_{C_\Pb(z)} (C_\Lb^\ci(v),z,v^{-1}u_z)$, which is equal to 
$C_\Gb^\ci(z) \times_{C_\Pb(z)} 
(\Sigb_z^\pr \times \Vb)_{\min,\Lb}^{C_\Gb^\ci(z)}$ (because $C_\Pb(z)$ 
is connected). Hence, it is contained 
in $\WC_z^\ci=W_{C_\Gb^\ci(z)}(\Lb,v)$, so it is equal 
to $H_{C_\Gb^\ci(z)}(\Lb,v,z)$ because this last variety is an open 
subset of $(\Xbt_{\min}^\pr)_\Lb^{C_\Gb^\ci(z)}$. 

Now, the action of $\WC_z^\ci$ on $(\Xbt_{\min}^\pr)_\Lb^{C_\Gb^\ci(z)}$ 
commutes with the translation by $z$. Hence $H_{C_\Gb^\ci(z)}(\Lb,v,z)=
H_{C_\Gb^\ci(z)}(\Lb,v,1)$.\fin

\bi

\rem The reader may be surprised by the fact that $C_\Gb(z)$ is not necessarily 
connected. However, 
he can check directly that the previous constructions ($\Ybh$, $\Xbt_\mini$, 
$W_\Gb(\Lb,v)$...) remains valid whenever $\Gb$ is not connected, provided 
that the \para $\Pb$ of $\Gb$ is connected.\finl

\bi

Proposition \ref{z} shows that it is necessary and sufficient 
to understand the stabilizer $H_\Gb(\Lb,v,1)$. 
However, we are only able to get a satisfying result 
whenever the centralizer of $u$ in $\Gb$ is contained in the \para $\Pb$.

\bi

\begin{prop}\label{agu}
If $C_\Gb(u) \incl \Pb$, then~:

\tete{1} $\pi_\mini^{\pr -1}(u)$ is the $A_\Lb(v)$-orbit of $\uti'$. In particular, 
$|\pi_\mini^{\pr -1}(u)|=|A_\Lb(v)|$.

\tete{2} $W_\Gb(\Lb,v)=A_\Lb(v) \rtimes H_\Gb(\Lb,v,1)$.
\end{prop}

\bi

\proof (1) follows immediately from Corollary \ref{agu 1}~: indeed, 
$\pi_\mini^{\pr -1}(u)=\fti^{-1}(\uti)$. 
By (1), $A_\Lb(v)$ acts freely and transitively on 
$\pi_\mini^{\pr -1}(u)$, so $W_\Gb(\Lb,v)$ acts transitively on 
$\pi_\mini^{\pr -1}(u)$. (2) follows from this remark and from 
\ref{H cap A}.\fin

\bi

\sub{Further investigations\label{01}} 
The group $C_\Gb^\ci(v) \cap \Lb=C_{C_\Gb^\ci(v)}(\Zb(\Lb)^\ci)$ 
is connected, because it is the centralizer of a torus in a connected 
group \cite[Corollary 11.12]{borel}. Therefore, we have the well-known equality 
\equat\label{cgo}
C_\Gb^\ci(v) \cap \Lb =C_\Lb^\ci(v).
\endequat
So the natural morphism $C_\Lb(v) \injto C_\Gb(v)$ induces an injective 
\mor 
\equat
A_\Lb(v) \injto A_\Gb(v).
\endequat
Let $W_\Gb^\ci(\Lb,v)=N_\Gb(\Lb,v) \cap C_\Gb^\ci(v)/C_\Lb^\ci(v)$. 
Since $C_\Gb^\ci(v) \cap \Lb=C_\Lb^\ci(v)$, we have $W_\Gb^\ci(\Lb,v) \cap 
A_\Lb(v)=1$. Moreover, $W_\Gb^\ci(\Lb,v)$ and $A_\Lb(v)$ are normal subgroups 
of $W_\Gb(\Lb,v)$, therefore $W_\Gb^\ci(\Lb,v) \times A_\Lb(v)$ is naturally 
a subgroup of $W_\Gb(\Lb,v)$. This discussion has the following immediate 
consequence~:

\bi

\begin{lem}\label{alv agv alv agv}
If $A_\Lb(v)=A_\Gb(v)$, then $W_\Gb(\Lb,v)=W_\Gb^\ci(\Lb,v) \times A_\Lb(v)$.
\end{lem}

\bi

\begin{coro}\label{mor fti}
Assume that $C_\Gb(u) \incl \Pb$, and that $A_\Lb(v)=A_\Gb(v)$. 
Then there exists a unique \mor of groups 
$\ph_{\Lb,v}^\Gb : W_\Gb^\ci(\Lb,v) \to A_\Lb(v)$ such that 
$$H_\Gb(\Lb,v,1)=\{(w,a) \in W_\Gb^\ci(\Lb,v) \times 
A_\Lb(v)~|~a=\ph_{\Lb,v}^{\Gb}(w)\}.$$
\end{coro}

\bi

\proof This follows from Proposition \ref{agu} (2) and from Lemma 
\ref{alv agv alv agv}.\fin

\bi

\begin{coro}\label{mor odd}
Assume that $C_\Gb(u) \incl \Pb$, that $A_\Lb(v)=A_\Gb(v)$, and that 
$|A_\Lb(v)|$ is odd. 
Then $H_\Gb(\Lb,v,1)=W_\Gb^\ci(\Lb,v)$.
\end{coro}

\bi

The morphism $\ph_{\Lb,v}^\Gb$ is the central object of this paper. 
In Part II, we will compute it explicitly whenever $v$ is a 
regular unipotent element under some restriction on $\Lb$. 

\bi

\sub{Separability} Let $\Cb^\et$ denote the separable closure of $\Cb$ in 
$\Lb/C_\Lb^\ci(v)$ (under the morphism $\Lb/C_\Lb^\ci(v) \to \Cb$, $l \mapsto 
lvl^{-1}$). Note that $\Cb^\et$ is smooth. 
The variety $\Cb^\et$ inherits from $\Lb/C_\Lb^\ci(v)$ 
the action of $\Lb$ by left translation, and the action of $A_\Lb(v)$ 
by right translation. Then we have a sequence of $\Lb \times A_\Lb(v)$-equivariant 
morphisms
$$\diagram
\Lb/C_\Lb^\ci(v) \rrto && \Cb^\et \rrto && \Cb.
\enddiagram$$
The first morphism is bijective and purely inseparable, the second one 
is a Galois \'etale covering with group $A_\Lb(v)$. We then define
$\Sigb^\et=\Zb(\Lb)^\ci \times \Cb^\et$, $\Sigb^\et_\reg=\Zb(\Lb)_\reg^\ci 
\times \Cb^\et$, $\Ybh^\et=\Gb \times \Sigb^\et_\reg$, and 
$\Ybt^\et=\Gb \times_\Lb \Sigb_\reg^\et$. We have a commutative diagram with 
cartesian squares
\equat\label{diagramme etale}
\diagram
\Sigb_\reg^\pr \ddto_{\DS{f_\reg^\ins}} && \Ybh' \llto_{\DS{\a'}} \rrto^{\DS{\b'}} 
\ddto_{\DS{\fha^\ins}} &&
 \Ybt' \ddto_{\DS{\fti^\ins}} \xto[4,2]^{\DS{\pi'}} & \\
 &&&&&&\\
\Sigb_\reg^\et \ddto_{\DS{f_\reg^\et}} && \Ybh^\et \llto_{\DS{\a^\et}} \rrto^{\DS{\b^\et}} 
\ddto_{\DS{\fha^\et}} &&
 \Ybt^\et \ddto_{\DS{\fti^\et}} \ddrrto_{\DS{\pi^\et}} & \\
&&&&&&\\
\Sigb_\reg && \Ybh \llto^{\DS{\a}} \rrto_{\DS{\b}} && \Ybt \rrto_{\DS{\pi}} && \Yb.
\enddiagram
\endequat
Here the maps $?^\et$ and $?^\ins$ are induced by the maps $?$ or $?'$. Moreover, 
all the morphisms $?^\et$ are Galois \'etale coverings, and all 
the morphisms $?^\ins$ are bijective purely inseparable morphisms.

By the same argument as in Remark \ref{xtilde 0}, the group $\Pb$ acts 
on the variety $\Sigb^\et \times \Vb$ and the quotient $\Xbt_0^\et=\Gb \times_\Pb 
(\Sigb^\et \times \Vb)$ exists~: it is the separable closure of $\Xbt_0$ 
in $\Xbt_0^\pr$. If we denote by $(\Sigb^\et \times \Vb)_\mini$ the inverse 
of $(\Sigb.\Vb)_\mini$ under the morphism $f^\et \times \Id_\Vb$, then 
$\Xbt_\mini^\et=\Gb \times_\Pb (\Sigb^\et \times \Vb)_\mini$ is the normalization 
of $\Xb_\mini$ in $\Ybt^\et$. So it inherits an action of $W_\Gb(\Lb,v)$ 
and the bijective purely inseparable morphism $\fti_\mini^\ins : \Xbt_\mini^\pr 
\to \Xbt_\mini^\et$ induced by $\fti_\mini$ 
is $W_\Gb(\Lb,v)$-equivariant. Moreover, the \mor $\fti_\mini^\et 
: \Xbt_\mini^\et \to \Xbt_\mini$ induced by $\fti_\mini$ 
is a Galois \'etale covering with group $A_\Lb(v)$. We summarize the notation 
in the next diagram.

$$
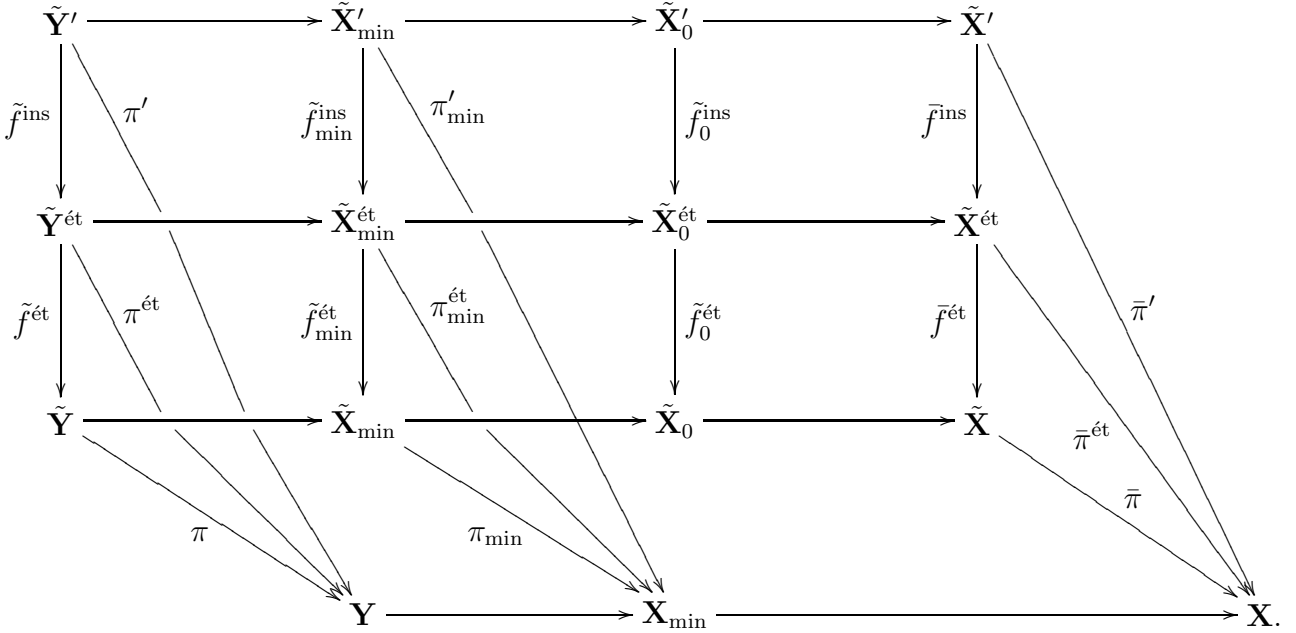
\diagram
\Ybt' \xto[0,3] \ddto_{\DS{\fti^\ins}}\xto'[2,1]^{\DS{\pi'}}'[4,2][6,3]&&& 
\Xbt_\mini^\pr \ddto_{\DS{\fti_\mini^\ins}}\xto'[2,1]^{\DS{\pi_\mini^\pr}}[6,3]
\xto[0,3] &&& \Xbt_0^\pr \ddto^{\DS{\fti_0^\ins}} \xto[0,3] &&& 
\Xbt' \ddto_{\DS{\fba^\ins}} \xto[6,3]^{\DS{\piba'}}\\
&&\\
\Ybt^\et \xto[0,3] \ddto_{\DS{\fti^\et}}\xto'[2,1]^{\DS{\pi^\et}}[4,3]&&& 
\Xbt_\mini^\et \ddto_{\DS{\fti_\mini^\et}}\xto'[2,1]^{\DS{\pi_\mini^\et}}[4,3]
\xto[0,3] &&& \Xbt_0^\et \ddto^{\DS{\fti_0^\et}} \xto[0,3] &&& 
\Xbt^\et \ddto_{\DS{\fba^\et}} 
\xto[4,3]_{\DS{\piba^\et}}&&& \\
&&\\
\Ybt \xto[0,3] \xto[2,3]_{\DS{\pi}} &&& \Xbt_\mini \xto[2,3]_{\DS{\pi_\mini}} \xto[0,3] 
&&& \Xbt_0  \xto[0,3] &&& \Xbt \xto[2,3]^{\DS{\piba}}\\
&&&\\
&&& \Yb \xto[0,3] &&& \Xb_\mini \xto[0,6] &&&&&&\Xb. \\
\enddiagram$$

\bi

\remark{stab etale} 
If $z \in \Zb(\Lb)^\ci$, we denote by $\uti_z^\et$ the image of $\uti_z^\pr \in 
\Xbt_\mini^\pr$ in $\Xbt_\mini^\et$ under the morphism $\fti_\mini^\ins$. 
Since $\fti_\mini^\ins$ is bijective and $W_\Gb(\Lb,v)$-equivariant, 
the stabilizer of $\uti_z^\et$ in $W_\Gb(\Lb,v)$ is equal to $H_\Gb(\Lb,v,z)$.\finl

\bi

\example{exemple etale} It may happen that the variety $\Cb^\et$ is different 
from $\Lb/C_\Lb^\ci(v)$, so that the construction above is not irrelevant. 
Of course, it only occurs in positive characteristic. The smallest example 
is given by the group $\Lb=\Gb=\Sb\Lb_2(\FM)$, whenever $p=2$ and
$$v=\matrice{1 & 1 \\ 0 & 1}.$$
Nevertheless, this is quite an unusual phenomenon. Indeed, 
if $\Gb=\Sb\Lb_n(\FM)$ and if $p$ does not divide $n$, then 
the morphism $\Lb/C_\Lb^\ci(v) \to \Cb$ is always \'etale. 
Also, if $\Gb$ is a quasisimple group of type different from $A$ and 
if $p$ is good for $\Gb$, then again the morphism $\Lb/C_\Lb^\ci(v) \to \Cb$ 
is always \'etale.\finl

\bi

\sec{Elementary properties of $\ph_{\Lb,v}^\Gb$\label{ele sec}}~

\med

As it will be shown in \SEC\ref{sec endo}, the knowledge of the morphism 
$\ph_{\Lb,v}^\Gb$ will be of fundamental 
importance in the description of the endomorphism algebra of an induced 
cuspidal character sheaf (\cf Corollary \ref{coro zeta phi}). That is the reason 
why we devote a section to gather the properties of this morphism. These 
properties may help to reduce the computations to small groups. 

\bi

\sub{Product of groups\label{produit}} Assume in this subsection, and 
only in this subsection, 
that $\Gb=\Gb_1 \times \Gb_2$. Let $\Lb=\Lb_1 \times \Lb_2$, $v=(v_1,v_2)$, 
$\Pb=\Pb_1 \times \Pb_2$. Then, for every $z =(z_1,z_2) \in \Zb(\Lb)^\ci$, we have
\equat
H_\Gb(\Lb,v,z)=H_{\Gb_1}(\Lb_1,v_1,z_1) \times H_{\Gb_2}(\Lb_2,v_2,z_2).
\endequat
Moreover, if $A_\Lb(v) =A_\Gb(v)$ and if $C_\Gb(u) \incl \Pb$, then 
\equat\label{produit H}
\ph_{\Lb,v}^\Gb = \ph_{\Lb_1,v_1}^{\Gb_1} \times \ph_{\Lb_2,v_2}^{\Gb_2}.
\endequat

\bi

\sub{Changing the group\label{changement}} Let $\Gb_1$ be a connected reductive 
algebraic group such that there exists a \mor of algebraic groups 
$\s : \Gb_1 \to \Gb$ such that

\sma

\begin{centerline}
{\begin{tabular}{l}
(1) {\it The kernel of $\s$ is central in $\Gb_1$,}\\
(2) {\it The image of $\s$ contains the derived group of $\Gb$.\hspace{3cm}}
\end{tabular}}\end{centerline}

\sma

\noindent We put $\Lb_1=\s^{-1}(\Lb)$ and $\Pb_1=\s^{-1}(\Pb)$, so that $\Lb_1$ 
is a \levi of the \para $\Pb_1$ of $\Gb_1$. Let $v_1$ denote the unique 
unipotent element of $\Gb_1$ such that $\s(v_1)=v$ and let $\Cb_1$ 
be the class of $v_1$ in $\Gb_1$. Note that $\s(\Cb_1)=\Cb$. 

\bi

\begin{lem}\label{G1 G}
$({\mathrm{a}})$ $\s(C_{\Lb_1}(v_1)) \Zb(\Gb)^\ci=C_\Lb(v)$.

\tete{b} The morphism $A_{\Lb_1}(v_1) \to A_\Lb(v)$ induced by $\s$ is surjective.

\tete{c} The morphism $W_{\Gb_1}(\Lb_1,\Sigb_1) \to W_\Gb(\Lb,\Sigb)$ induced by $\s$ 
is an isomorphism.

\tete{d} The morphism $W_{\Gb_1}(\Lb_1,v_1) \to W_\Gb(\Lb,v)$ induced by $\s$ 
is surjective.

\tete{e} The morphism $W_{\Gb_1}^\ci(\Lb_1,v_1) \to W_\Gb^\ci(\Lb,v)$ induced by $\s$ 
is an isomorphism.
\end{lem}

\bi

\proof First, note that (b) is an immediate consequence of (a). 
Let $l \in C_\Lb(v)$. Then there exists $l_1$ in $\Lb_1$ and $z \in \Zb(\Gb)^\ci$ such 
that $\s(l_1)=lz$. So there exists $x \in \Ker \s$ such that $l_1v_1l_1^{-1}=xv_1$. 
But $l_1v_1l_1^{-1}$ is unipotent and $\Ker \s$ consists of central 
semisimple elements of $\Lb_1$, so $x=1$, that is 
$l_1$ centralizes $v_1$. Hence (a) follows.

(c), (d) and (e) follow by similar arguments (note that 
$\s(C_{\Gb_1}^\ci(v_1)).\Zb(\Gb)^\ci=C_\Gb^\ci(v)$).\fin 

\bi

We will denote by a subscript $?_1$ the object associated to the datum $(\Lb_1,v_1)$ 
and defined in the same way as the object $?$ in $\Gb$ (for instance, $\Sigb_1=
\Zb(\Lb_1)^\ci .\Cb_1$, $\Xbt_1$, $\Ybt_1^\pr$...).

The reader may check that $\s((\Sigb_1.\Vb_1)_\mini) \incl (\Sigb.\Vb)_\mini$, so 
that $\s$ induces a morphism $\Xbt_{1,\min} \to \Xbt_\mini$. This morphism 
is $W_{\Gb_1}(\Lb_1,\Sigb_1)$-equivariant as it can be checked by restriction to 
$\Ybt_1$ (here, $W_{\Gb_1}(\Lb_1,v_1)$ and $W_\Gb(\Lb,\Sigb)$ are identified 
via the morphism $\s$ by Lemma \ref{G1 G} (c)). 
Similarly, $\s$ induces a $W_{\Gb_1}(\Lb_1,v_1)$-equivariant 
morphism $\Xbt_{1,\min}^\pr \to \Xbt_\mini^\pr$ (here, $W_{\Gb_1}(\Lb_1,v_1)$ 
acts on $\Xbt_\mini^\pr$ via the surjective morphism $W_{\Gb_1}(\Lb_1,v_1)
\to W_\Gb(\Lb,v)$). 

\bi

\begin{prop}\label{G1 G H}
If $z_1 \in \Zb(\Lb_1)^\ci$, then~:

\tete{a} $H_{\Gb_1}(\Lb_1,\Sigb_1,z_1)=
H_\Gb(\Lb,\Sigb,\s(z_1))$. 

\tete{b} $\s$ induces an isomorphism 
$H_{\Gb_1}(\Lb_1,v_1,z_1)\simeq H_\Gb(\Lb,v,\s(z_1))$.

Moreover,

\tete{c} If $C_{\Gb_1}(u_1) \incl \Pb_1$ and if $A_{\Lb_1}(v_1)=A_{\Gb_1}(v_1)$, 
then $C_\Gb(u) \incl \Pb$ and $A_\Lb(v)=A_\Gb(v)$. In this case, 
the diagram
$$\diagram
W_{\Gb_1}^\ci(\Lb_1,v_1) \rrto_{\sim}^{\DS{\s}} \ddto_{\DS{\ph_{\Lb_1,v_1}^{\Gb_1}}}
&& W_\Gb^\ci(\Lb,v) \ddto^{\DS{\ph_{\Lb,v}^\Gb}} \\
&&\\
A_{\Lb_1}(v_1) \rrto^{\DS{\s}} && A_\Lb(v) 
\enddiagram$$
is commutative.
\end{prop}

\bi

\proof It is clear that $H_{\Gb_1}(\Lb_1,\Sigb_1,z_1) \incl 
H_\Gb(\Lb,\Sigb,\s(z_1))$. To prove the reverse inclusion, we use Proposition 
\ref{stabilisateur z} 
to reduce the problem to the case where $z_1=1$. Then if an element 
$w \in W_\Gb(\Lb,\Sigb)$ stabilizes $1 *_\Pb u$, this proves that there exists 
$z_1 \in \Ker \s$ such that $(1 *_{\Pb_1} u_1).w=1*_{\Pb_1} z_1u_1$. But, by 
\ref{ZZZ}, $z_1=1$. This proves (a).

To prove (b), we must notice that 
$\s(H_{\Gb_1}(\Lb_1,v_1,z_1) )\incl H_\Gb(\Lb,v,\s(z_1))$, so $\s$ induces 
a \mor $H_{\Gb_1}(\Lb_1,v_1,z_1) \longto  H_\Gb(\Lb,v,\s(z_1))$. This morphism 
is injective by \ref{H cap A}. It is surjective by the same argument as the 
one used in (a). The proof of (b) is complete.

The first assertion of (c) follows from Lemma \ref{G1 G}, (a), while 
the second follows from (b).\fin

\bi

\sub{Parabolic restriction\label{para sub red}} In this subsection, we show that 
the above constructions are compatible with ``restriction'' to parabolic subgroups 
of $\Gb$. We need some notation. Let $\Pb'$ denote a \para of $\Gb$ containing 
$\Pb$ and let $\Lb'$ denote the unique \levi of $\Pb'$ containing $\Lb$. 
Let $\Vb'$ denote the unipotent radical of $\Pb'$, and let $v' =
\pi_{\Lb'}(u)$. Then $v' \in v(\Vb \cap \Lb')$. We start by some 
elementary properties.

\bi

\begin{prop}\label{fourre-tout}
$(1)$ $v' \in \Cb^{\Lb'}$.

\tete{2} If $A_\Lb(v) =A_\Gb(v)$, then $A_\Lb(v)=A_{\Lb'}(v)$. 

\tete{3} If $C_\Gb(u) \incl \Pb$, then $C_{\Lb'}(v') \incl \Pb \cap \Lb'$.
\end{prop}

\proof By Lemma \ref{dim p}, we have $\dim C_\Pb(u) \ge \dim C_{\Pb \cap \Lb'}(v') \ge 
\dim C_\Lb(v)$. But 
$$\dim C_\Pb(u)=\dim C_\Lb(v)$$ 
because $u \in \Cb^\Gb$. 
So, $\dim C_{\Pb \cap \Lb'}(v') = \dim C_\Lb(v)$. This proves 
that $\dim~(v')_{\Pb \cap \Lb'} = \dim \Cb + \dim \Vb \cap \Lb'$. 
So $(v')_{\Pb \cap \Lb'}$ is dense in $\Cb.(\Vb \cap \Lb')$, that is $v' \in \Cb^{\Lb'}$. 
Hence (1) is proved. 

(2) follows from \ref{cgo} applied to the \levi $\Lb'$. Let us now prove 
(3). Let $m \in C_\Lb'(v')$. We only need to prove that $m \in \Pb$. 
But $\lexp{m}{u} \in v'.\Vb' \cap \Cb^\Gb \incl v\Vb \cap \Cb^\Gb$. So, 
by \cite[Proposition II.3.2 (d)]{spaltenstein}, there exists $x \in \Pb$ 
such that $\lexp{m}{u}=\lexp{x}{u}$. So $x^{-1}m \in C_\Gb(u)$. 
But $C_\Gb(u) \incl \Pb$ by hypothesis, so $m \in \Pb$.\fin

\bi

\begin{prop}\label{prop compa para}
If $C_\Gb(u) \incl \Pb$, then $H_\Gb(\Lb,v,1) \cap W_{\Lb'}(\Lb,v)=H_{\Lb'}(\Lb,v,1)$.
\end{prop}

\bi

\proof By Proposition \ref{agu} (2), the subgroups 
$H_\Gb(\Lb,v,1) \cap W_{\Lb'}(\Lb,v)$ and $H_{\Lb'}(\Lb,v,1)$ 
have the same index in $W_{\Lb'}(\Lb,v)$ (this index is equal to $|A_\Lb(v)|$). 
Consequently, it is sufficient to prove that 
$$H_\Gb(\Lb,v,1) \cap W_{\Lb'}(\Lb,v) \incl H_{\Lb'}(\Lb,v,1).\leqno{(\#)}$$
Let 
$$\Fbt'=\Pb' \times_\Pb (\Sigb^\pr \times \Vb)_{\min,\Lb}^\Gb.$$
It is an irreducible closed subvariety of $\Xbt'$, and it is stable 
under the action of $W_{\Lb'}(\Lb,v)$ (indeed, the open subset $\Obt'=\Pb' \times_\Lb 
\Sigb_\reg^\pr$ is obviously $W_{\Lb'}(\Lb,v)$-stable). 

By Lemma \ref{dim p} (2), the projection $\pi_{\Lb'} : \Pb' \to \Lb'$ sends an element of 
$(\Sigb.\Vb)_{\mini,\Lb}^\Gb$ to an element of 
$(\Sigb.(\Vb \cap \Lb'))_{\min,\Lb}^{\Lb'}$. 
So, it induces a map $\g : \Fbt' \to (\Xbt_\mini^\pr)_\Lb^{\Lb'}$. 
Moreover, the diagram
$$\diagram
\Obt' \rrto \ddto &&\Fbt' \ddto^{\DS{\g}} \\
&&\\
(\Ybt')_\Lb^{\Lb'} \rrto && (\Xbt_\mini^\pr)_\Lb^{\Lb'}
\enddiagram$$
is commutative. The first vertical map is $W_{\Lb'}(\Lb,v)$-equivariant, so, by density, 
the second vertical map is also $W_{\Lb'}(\Lb,v)$ equivariant. This 
proves $(\#)$.\fin

\bi

\begin{coro}\label{coro compa para}
If $C_\Gb(u) \incl \Pb$ and if $A_\Lb(v)=A_\Gb(v)$, then
$$\ph_{\Lb,v}^{\Lb'}=\Res_{W_{\Lb'}^\ci(\Lb,v)}^{W_\Gb^\ci(\Lb,v)} \ph_{\Lb,v}^\Gb.$$
\end{coro}

\bi

\proof By Proposition \ref{fourre-tout}, (2) and (3), 
$\ph_{\Lb,v}^{\Lb'}$ is well-defined. 
So Corollary \ref{coro compa para} follows from 
Proposition \ref{prop compa para}.\fin

\bi

\remark{qwerty} If $\Gb'$ is a connected reductive subgroup of $\Gb$ containing 
$\Lb$, then it may happen that 
$$\ph_{\Lb,v}^{\Gb'}\not=
\Res_{W_{\Gb'}^\ci(\Lb,v)}^{W_\Gb^\ci(\Lb,v)} \ph_{\Lb,v}^\Gb.$$
An example is provided in Part II of this paper.\finl

\bi


\bi

\sec{Endomorphism algebra of induced cuspidal character sheaves\label{sec endo}}~

\med

\noindent{\bf Hypothesis and notation :} {\it From now on, and until the end 
of this first part, we assume that $\Cb$ supports an irreducible cuspidal 
\cite[Definition 2.4]{luicc} local system $\EC$. 
To this local system is associated an irreducible character $\z$ of 
$A_\Lb(v)$, via the Galois \'etale covering $\Cb^\et \to \Cb$. 
Let $\FC = \qlb \boxtimes \EC$ ($\FC$ is a local system on $\Sigb$) and let 
$\FC_\reg$ denote the restriction of $\FC$ to $\Sigb_\reg$.
Let $K$ be the perverse sheaf on $\Gb$ obtained from the triple $(\Lb,v,\z)$ 
by parabolic induction \cite[4.1.1]{lucs} and let $\AC$ denote its endomorphism 
algebra.}

\bi

In \cite[Theorem 9.2]{luicc}, Lusztig constructed an isomorphism 
$\Th : \qlb W_\Gb(\Lb) \to \AC$. This \iso is very 
convenient for computing the generalized Springer correspondence. 
On the other hand, Lusztig's construction is canonical but not explicit. 
The principal aim of this paper is to construct 
an explicit \iso $\Th' : \qlb W_\Gb(\Lb) \to \AC$ by another 
method. It turns out that this \iso differs from Lusztig's one 
by a linear character $\g$ of $W_\Gb(\Lb)$. The knowledge of $\g$ 
would allow us to combine the advantages of both 
\isos $\Th$ and $\Th'$. However, we are not able to determine 
it explicitly in general, although we can relate it 
with the morphism $\ph_{\Lb,v}^\Gb$ defined in the previous section. 
Note that $\g=1$ whenever $\Lb=\Tb$. 

In this section we recall some well-known facts about cuspidal 
local systems, parabolic induction and endomorphism algebra. 
Most of these results may be found in \cite{luicc} or \cite{lucs}. 
However, Theorem \ref{alv}, which is of fundamental step 
for constructing the \iso $\Th'$ defined above, 
was proved in full generality in \cite{Bonnafe}. 
The \iso $\Th'$ will be constructed in \SEC\ref{section endo}. 
We will also prove in \SEC\ref{section endo} the existence of $\g$ 
and its relation with $\ph_{\Lb,v}^\Gb$. 
We will provide in \SEC\ref{section 1 gamma} some properties of $\g$ which 
allows us to reduce its computation in the case where $\Gb$ is semisimple, 
simply-connected, quasi-simple, and $\Pb$ is a maximal \para of $\Gb$. 
In Part II, we will use the knowledge of the morphism $\ph_{\Lb,v}^\Gb$ 
for $v$ regular to determine explicitly the linear character 
$\g$. However, we will need to assume that $p$ is good for $\Lb$. 
As an application, 
we will get some precision on Digne, Lehrer and Michel's theorem of 
Lusztig restriction of Gel'fand-Graev characters \cite[Theorem 3.7]{DLM2}.

\bi

\sub{Parabolic induction\label{subsec para}} For the convenience of the reader, 
we reproduce here the diagram \ref{diagramme central}.
$$\diagram
\Sigb_\reg \ddto && \Ybh \llto_{\DS{\a}} 
\rrto^{\DS{\b}} \ddto && \Ybt \rrto^{\DS{\pi}} \ddto && \Yb \ddto\\
&&&&&&\\
\overline{\Sigb} && \Xbh \llto_{\DS{\alpba}} \rrto^{\DS{\betba}} && \Xbt 
\rrto^{\DS{\piba}}  && \Xb.
\enddiagram$$
We define $\FCh_\reg = \a^* \FC_\reg$~: it is a local system on $\Ybh$. 
Moreover, since $\EC$ is $\Lb$-equivariant, there exists a local system $\FCt_\reg$ 
on $\Ybt$ such that $\b^* \FCt_\reg \simeq \FCh_\reg$.
By \cite[3.2]{luicc}, 
the \mor $\pi : \Ybt \to \Yb$ is a Galois covering with Galois group $W_\Gb(\Lb)
=N_\Gb(\Lb)/\Lb$, 
so $\pi_*\FCt_\reg=\pi_!\FCt_\reg$ (because $\pi$ is finite hence proper) is a 
local system on $\Yb$. We denote by $K$ the following perverse sheaf on $\Gb$~:
\equat
K=IC(\overline{\Yb},\KC)[\dim \Yb].
\endequat
where $\pi_* \FCt_\reg=\KC$. Recall that $\overline{\Yb}=\Xb$, so that $\dim \Yb
 = \dim \Xb$.

We shall give, following \cite[\SEC 4]{luicc}, 
an alternative description of the perverse sheaf $K$. 
Let $A$ be the following perverse sheaf on $\Lb$~:
$$A=IC(\overline{\Sigb},\FC)[\dim \Sigb].$$
Note that $\overline{\Sigb}=\Zb(\Lb)^\ci \overline{\Cb}$ so $A=\qlb[\dim \Zb(\Lb)^\ci] 
\boxtimes IC(\overline{\Cb},\EC)[\dim \Cb]$. 
Since $A$ is $\Lb$-equivariant, 
there exists a perverse sheaf $\Kti$ on $\Xbt$ such that 
$$\alpba^*A[\dim \Gb + 
\dim \Vb] = \betba^*\Kti[\dim \Pb].$$
The perverse sheaf $\Kti$ is in fact 
equal to $IC(\Xbt,\FCt_\reg)[\dim \Xbt]$. By \cite[Proposition 4.5]{luicc}, we have
\equat\label{K}
K={\mathrm{R}}\piba_*\Kti.
\endequat

\bi

\sub{Recollection\label{sub reco}} The fact that $\Cb$ admits a cuspidal 
local system has a lot of consequences. We gather some of them in the next 
theorem.

\bi

\begin{theo}[{\bf Lusztig}]\label{rappellusztig}
\tete{a} $v$ is a {\bf distinguished} unipotent element of $\Lb$ $($that is, 
$v$ is not contained in a \levi of a proper \para of $\Lb)$.

\tete{b} $N_\Gb(\Lb)$ stabilizes $\Cb$ and $\EC$.
\end{theo}

\bi

\proof \cf \cite[Proposition 2.8]{luicc} for (a) and \cite[Theorem 9.2]{luicc}  
for (b) and (c).\fin

\bi

The next theorem has been proved by Lusztig provided that $p$ is large enough by using 
the classification of cuspidal pairs \cite{luhecke}. In 
\cite[Corollary to Proposition 1.1]{Bonnafe}, the author gave a proof 
in the general case without using the classification.

\bi

\begin{theo}\label{alv}
The injective \mor $C_\Lb(v) \to C_\Gb(v)$ induces an \iso of 
finite groups $A_\Lb(v) 
\to A_\Gb(v)$.
\end{theo}

\bi

By Theorem \ref{alv} and Lemma \ref{alv agv alv agv}, we have 
\equat\label{W produit}
W_\Gb(\Lb,v)=A_\Lb(v) \times W_\Gb^\ci(\Lb,v),
\endequat
and, by Theorem \ref{rappellusztig} (b), we have 
\equat\label{W}
W_\Gb^\ci(\Lb,v) \simeq W_\Gb(\Lb)=N_\Gb(\Lb)/\Lb.
\endequat

\bi

\sub{Lusztig's description of $\AC$\label{sub endo lusztig}} For each $w$ in
$W_\Gb^\ci(\Lb,v)$, we choose a representative $\wdo$ of $w$ in $N_\Gb(\Lb) 
\cap C_\Gb^\ci(v)$. By Theorem \ref{rappellusztig} (b), the local 
systems $\FC$ and $(\INT \dot{w})^* \FC$ are isomorphic. Let $\th_w$ denote an 
\iso of $\Lb$-equivariant local systems $\FC \to (\INT \dot{w})^* \FC$.
Then $\th_w$ induces an \iso $\thet_w : \FCt_\reg \to \g_w^* \FCt_\reg$ where 
$\g_w : \Ybt \to \Ybt$, $(g,x\Lb) \mapsto (g,x\wdo^{-1}\Lb)$ 
(\cf \cite[proof of Proposition 3.5]{luicc}). 
But $\pi_* \g_w^* = \pi_*$. Hence $\pi_*\thet_w$ is an \auto of $\KC$. 
By applying the functor $IC(\Xb,.)[\dim \Yb]$, $\pi_* \thet_w$ induces 
an \auto $\Th_w$ of $K$. The \auto $\Th_w$, as well as $\th_w$, is 
defined up to multiplication by an element of $\qlb^\times$. 
By \cite[Proposition 3.5]{luicc}, $(\Th_w)_{w \in W_\Gb^\ci(\Lb,v)}$ is a basis 
of the endomorphism algebra $\AC$ of $K$~; 
moreover, by \cite[Remark 3.6]{luicc}, 
there exists a family of scalar $(a_{w,w'})_{w,w' \in W_\Gb^\ci(\Lb,v)}$ 
of elements of $\qlb^\times$ 
such that $\Th_w \Th_{w'}= a_{w,w'} \Th_{ww'}$ for all $w$ and $w'$ in 
$W_\Gb^\ci(\Lb,v)$. 
Lusztig proved that it is possible to choose in a 
canonical way the family $(\th_w)_{w \in W_\Gb^\ci(\Lb,v)}$ 
such that $\Th_w \Th_{w'}= \Th_{ww'}$ for all $w$ and $w'$ in 
$W_\Gb^\ci(\Lb,v)$. The next theorem \cite[Theorem 9.2]{luicc} 
explains his construction. 

\bi

\begin{theo}[{\bf Lusztig}]\label{lusztig theta}
There exists a unique family of \isos of local systems 
$(\th_w : \FC \to (\INT \dot{w})^* \FC)_{w \in W_\Gb^\ci(\Lb,v)}$ 
such that the following condition holds~:  
for each $w \in W_\Gb^\ci(\Lb,v)$, $\Th_w$ acts trivially 
on $\HC^{- \dim \Yb}_u K$, where $u$ is any element of $\Cb^\Gb$. 
\end{theo}

\bi

In the previous theorem, the uniqueness of the family 
$(\th_w)_{w \in W_\Gb^\ci(\Lb,v)}$ follows from the fact that 
$\HC^{- \dim \Yb}_u K \not= 0$ for each $u \in \Cb^\Gb$. 
As a consequence, one gets that the linear mapping 
$$\fonction{\Th}{\qlb W_\Gb^\ci(\Lb,v)}{\AC=\End_{\MC\Gb}(K)}{w}{\Th_w}$$
is an \iso of algebras. 
If $\ch$ is an \irr \car of $W_\Gb^\ci(\Lb,v)$, we denote by $K_\ch$ 
an irreducible component of $K$ associated to $\ch$ via the 
\iso $\Th$. 

\bi

\begin{coro}[{\bf Lusztig}]\label{lusztig coro}
For each $u \in \Cb^\Gb$, we have~:

\tete{a} $\HC^{-\dim \Yb}_u K_1 = \HC^{-\dim \Yb}_u K$. 

\tete{b} $\HC^{-\dim \Yb}_u K_\ch = 0$ for every non-trivial \irr \car $\ch$ of
$W_\Gb^\ci(\Lb,v)$.
\end{coro}

\bi

\sub{Restriction to the open subset $\Xb_\mini$} 
The restriction $\KCt_0$ of $\Kti[-\dim \Xb]$ to $\Xbt_0$ 
is a local system \cite[4.4]{luicc}, 
that is, a complex concentrated in degree $0$, whose $0$th term 
is a local system. In fact, $\KCt_0$ is the local system on $\Xbt_0$ 
associated to the Galois \'etale covering $\Xbt_0^\et \to \Xbt_0$ and to 
the character $\z$ of $A_\Lb(v)$. Therefore, the restriction $\KCt_\mini$ 
of $K[-\dim \Xb]$ to $\Xbt_\mini$ is a local system. More precisely, 
it is the local system associated to the Galois \'etale covering 
$\Xbt_\mini^\et \to \Xbt_\mini$ and to the character $\z$.
Let $K_{\min}$ denote the restriction of $K[- \dim \Xb]$ to $\Xb_{\min}$. 
We have the following result.

\bi

\begin{prop}\label{k min}
We have $K_{\min}=\pi_{\min,*} \Kti_{\min}$. So, $K_{\min}$ is a 
constructible sheaf, that is a complex concentrated in degree $0$.
\end{prop}

\bi

\proof Since $\pi_{\min}$ is finite, the functor $\pi_{\min,*}$ is exact. 
The proposition follows from this remark and the Proper Base Change Theorem.\fin

\bi

\sec{Another isomorphism between $\AC$ and 
$\qlb W_\Gb^\ci(\Lb,v)$\label{section endo}} 

\bi

The aim of this section is to construct an explicit \iso 
$\Th'$ between the endomorphism algebra $\AC$ and 
the group algebra $\qlb W_\Gb^\ci(\Lb,v)$. Our strategy is 
the following. First, note that the endomorphism algebra $\AC$ of $K$ is 
canonically isomorphic to the endomorphism algebra of the local 
system $\KC$ on $\Yb$. To this local system is associated 
a representation of the fundamental group $\pi_1(\Yb,y)$ of $\Yb$ 
(here, $y$ is any point of $\Yb$). This representation and its endomorphism 
algebra are easy to describe (\cf \ref{iso psi}). 

\bi

\sub{Representations of the fundamental group} 
Let $V_\z$ denote an \irr left $\qlb A_\Lb(v)$-module affording the character $\z$. 
We may, and we will, assume that
$$\FC=(f^\et)_* \qlb \otimes_{\qlb A_\Lb(v)} V_\z,$$
$$\FC_\reg=(f_\reg^\et)_* \qlb \otimes_{\qlb A_\Lb(v)} V_\z,$$
$$\FCh_\reg=(\fha_\reg^\et)_* \qlb \otimes_{\qlb A_\Lb(v)} V_\z,$$
$$\FCt_\reg=(\fti_\reg^\et)_* \qlb \otimes_{\qlb A_\Lb(v)} V_\z.\leqno{\mathrm{and}}$$
Here, $V_\z$ is identified with the constant sheaf with values in $V_\z$. 
From the third equality, we deduce that
$$\KC=\pi_* \FCt_\reg=\pi_*^\et \qlb \otimes_{\qlb W_\Gb(\Lb)} 
\Ind_{A_\Lb(v)}^{W_\Gb(\Lb,v)} V_\z.$$
Therefore, the \endo algebra of $\KC$ is canonically isomorphic to the endomorphism 
algebra of the $\qlb W_\Gb(\Lb,v)$-module $\Ind_{A_\Lb(v)}^{W_\Gb(\Lb,v)} V_\z$. 
But, by \ref{W produit}, this endomorphism algebra is canonically isomorphic 
to $\qlb W_\Gb^\ci(\Lb,v)$. 
Since the functor $IC(\overline{\Yb}, .)[\dim \Yb]$ 
is fully faithful, it induces an \iso 

\equat\label{iso psi}
\Th' : \qlb W_\Gb^\ci(\Lb,v) \longto \AC.
\endequat

\sma

This \iso may be constructed in another way. 
The action of an element $\wdo \in N_\Gb(\Lb) \cap C_\Gb^\ci(v)$ on 
$\Sigb_\reg^\et$, $\Ybh^\et$, and $\Ybt^\et$ commutes with the action of 
$A_\Lb(v)$. Therefore, there exists an \iso $\th_w^\pr : 
\FC \to (\INT \dot{w})^* \FC$ (\resp $\t_w^\pr : 
\FCt_\reg \to (\INT \dot{w})^* \FCt_\reg$) which induces the identity on the stalks 
at $zv$ (\resp $1 *_\Lb zv$) for every $z \in \Zb(\Lb)^\ci$ (\resp for every 
$z \in \Zb(\Lb)^\ci_\reg$). Then 
\equat\label{tau theta}
\b^*\t_w^\pr = \a^* \th_w^\pr|_{\Sigb_\reg}.
\endequat
Now, let $\Th_w^\pr =IC(\overline{\Yb},\pi_* \t_w^\pr)[\dim \Yb] : K 
\longmapright{\sim} K$. By \ref{tau theta}, 
there exists an element $\g_w \in \qlb^\times$ such that
\equat\label{gamma g}
\Th(w)=\g_w \Th_w^\pr.
\endequat
By looking at the action on the stalk at $zv \in \Yb$, one can immediately get 
the following result.

\bi 

\begin{prop}\label{wouaw}
With the above notation, we have $\Th^\pr_w=\Th'(w)$ for every 
$w \in W_\Gb^\ci(\Lb,v)$.
\end{prop}

\bi

\begin{coro}\label{epsilon}
There exists 
a linear character $\g_{\Lb,v,\z}^\Gb$ of the Weyl group $W_\Gb^\ci(\Lb,v)$ such that
$$\Th^\pr(w)=\g_{\Lb,v,\z}^\Gb(w) \Th(w)$$
for every $w \in W_\Gb^\ci(\Lb,v)$.
\end{coro}

\bi

\proof This follows from Theorem \ref{lusztig theta}, from \ref{gamma g}, 
and from Proposition \ref{wouaw}.\fin

\bi
 
If $\ch$ is an \irr \car of $W_\Gb^\ci(\Lb,v)$, we 
denote by $K_\ch^\pr$ an irreducible component of $K$ associated to $\ch$ 
via the \iso $\Th'$. By Corollary \ref{epsilon}, we have
\equat\label{gamma idiot}
K_\ch^\pr=K_{\g_{\Lb,v,\z}^\Gb \ch}.
\endequat

\bi

\sub{Links between $\g_{\Lb,v,\z}^\Gb$ and $\ph_{\Lb,v}^\Gb$} 
The following proposition is an immediate consequence of 
\ref{gamma idiot} and Corollary \ref{lusztig coro}~:

\bi

\begin{prop}\label{carac gamma}
The linear character $\g_{\Lb,v,\z}^\Gb$ of $W^\ci_\Gb(\Lb,v)$ is the 
unique \irr \car $\g$ of $W_\Gb^\ci(\Lb,v)$ 
satisfying $\HC_u^{-\dim \Yb} K_\g^\pr \not= 0$ for some (or any) $u \in \Cb^\Gb$.
\end{prop}

\bi

If $\ch$ is an irreducible \car of $W_\Gb^\ci(\Lb,v)$, we denote by 
$K_{\mini,\ch}$ (\resp $K_{\mini,\ch}^\pr$) the irreducible component 
of $K_{\mini}$ associated to $\ch$ via the \iso $\Th$ (\resp $\Th'$). 
Now, let $V_\ch$ be an irreducible $\qlb W_\Gb^\ci(\Lb,v)$-module 
affording $\ch$ 
as character. Then, since $W_\Gb(\Lb,v)$ acts on $\Xbt_\mini^\et$, 
it also acts on the constructible sheaf $(\pi_\mini^\et)_* \qlb$ and we have 
by construction
\equat\label{equation K}
K_{\min,\ch}^\pr=(\piba^\et_* \qlb)^{\mathrm{opp}} \otimes_{\qlb W_\Gb(\Lb,v)} 
(V_\ch \otimes V_\z).
\endequat
Since $u \in \Xb_\mini$, the stalk of $\HC_u^{-\dim \Yb} K_\ch^\pr$ at $u$ may easily 
be deduced from \ref{equation K}~: we have
$$\dim_{\qlb} \HC_u^{-\dim \Yb} K_\ch^\pr= \dim (K_{\min,\ch}^\pr)_u = 
\langle 
\Res_{H_\Gb(\Lb,v,1)}^{W_\Gb(\Lb,v)} (\ch \otimes \z), \boldsymbol{1}_{H_\Gb(\Lb,v,1)} 
\rangle$$
We deduce immediately from this 
and from Proposition \ref{carac gamma} the following result.

\bi

\begin{prop}\label{zeta phi}
The linear character $\g_{\Lb,v,\z}^\Gb$ is the unique irreducible 
character $\g$ of $W_\Gb^\ci(\Lb,v)$ such that $\langle 
\Res_{H_\Gb(\Lb,v,1)}^{W_\Gb(\Lb,v)} (\g \otimes \z), \boldsymbol{1}_{H_\Gb(\Lb,v,1)} 
\rangle \not= 0$. 
\end{prop}

\bi

\begin{coro}\label{coro zeta phi}
If $C_\Gb(u) \incl \Pb$, then 
$$\g_{\Lb,v,\z}^\Gb={1 \over \z(1)} \z \ci \ph_{\Lb,v}^\Gb.$$
\end{coro}

\bi

\proof Since $A_\Lb(v)=A_\Gb(v)$, the morphism $\ph_{\Lb,v}^\Gb : W_\Gb^\ci(\Lb,v) 
\to A_\Lb(v)$ is well-defined. Hence, the corollary follows from 
Proposition \ref{zeta phi} and Corollary \ref{mor fti}.\fin

\bi

\begin{coro}\label{odd}
If $G_\Gb(u) \incl \Pb$, and if $|A_\Lb(v)|$ is odd, then $\g_{\Lb,v,\z}^\Gb=1$.
\end{coro}

\bi

\begin{coro}\label{tore}
$\g_{\Tb,1,\qlb}^\Gb=1$.
\end{coro}

\bi

\proof Indeed, if $\Pb=\Bb$, then $u$ is a regular unipotent element and $v=1$. 
Hence $C_\Gb(u) \incl \Bb$ so Corollary \ref{coro zeta phi} 
applies. But $A_\Lb(v)=1$, so $\g_{\Tb,1,\qlb}^\Gb=1$.\fin

\bi

\example{regular example} Assume in this example that $v$ is a {\it regular} 
unipotent element of $\Lb$. In this case, $u$ is a regular unipotent element 
of $\Gb$, so $C_\Gb(u) \incl \Pb$. Moreover, $A_\Lb(v)$ is abelian \cite{springer}. 
So $\z$ is a linear character and Corollary \ref{coro zeta phi} applies. 
We get
$$\g_{\Lb,v,\z}^\Gb=\z \ci \ph_{\Lb,v}^\Gb.$$
This case will be studied in full details in Part II.\finl

\bi

%
%
\sec{Elementary properties of the character $\g_{\Lb,v,\z}^\Gb$\label{section 1 gamma}}~

\med

\sub{Product of groups} We assume in this subsection that $\Gb=\Gb_1 \times \Gb_2$ 
where $\Gb_1$ and $\Gb_2$ are reductive groups. Let $\Lb=\Lb_1 \times 
\Lb_2$, $v=(v_1,v_2)$, $\Cb=\Cb_1 \times \Cb_2$ and $\z=\z_1 \otimes \z_2$. 
Then it is clear that 
\equat
W_\Gb^\ci(\Lb,v) = W_{\Gb_1}^\ci(\Lb_1,v_1) \times W_{\Gb_2}^\ci(\Lb_2,v_2)
\endequat
and that 
\equat\label{gamma produit}
\g_{\Lb,v,\z}^\Gb = \g_{\Lb_1,v_1,\z_1}^{\Gb_1} \otimes 
\g_{\Lb_2,v_2,\z_2}^{\Gb_2}.
\endequat

\bi

\sub{Changing the group\label{subsection change}} We use here the notation introduced 
in \SEC\ref{changement}. In particular, we still denote by a subscript $?_1$ for the object in 
$\Gb_1$ corresponding to $?$ in 
$\Gb$ (e.g. $\Lb_1$, $\FC_1$, $\FC_{\reg,1}$, $\Ybh_1$, $K_1$, $\AC_1$, $\Th_1$\dots). 
We have $\s^{-1}(\Yb)=\Yb_1$. 
Note that the groups $W_{\Gb_1}^\ci(\Lb_1,v_1)$ and $W_\Gb^\ci(\Lb,v)$ are 
isomorphic via $\s$ by Lemma \ref{G1 G} (e). 

\bi

\begin{lem}
The local system $\s^*\EC_1$ on $\Cb_1$ is cuspidal. It is associated 
to the irreducible character $\z_1$ of $A_{\Lb_1}(v_1)$ obtained from $\z$ by composing 
with the surjective \mor $A_{\Lb_1}(v_1) \to A_\Lb(v)$ $($\cf Lemma 
\ref{G1 G} ${\mathrm{(b))}}$.
\end{lem}

\bi

\proof 
This is immediate from the alternative definition of a cuspidal local system given 
in terms of permutation representations \cite[introduction]{luicc}.\fin

\bi

\begin{prop}\label{K1}
\tete{a} The restriction of $\s^*(K)$ to $\Yb_1$ is isomorphic to $K_1$.

\tete{b} $\s$ induces an isomorphism $\dot{\s} : \AC_1 \simeq \AC$.

\tete{c} The diagrams
$$\diagram
\qlb W_{\Gb_1}^\ci(\Lb_1,v_1) \rrto^{\qquad\DS{\Th_1}} \ddto_{\DS{\s}}&& \AC_1
\ddto^{\DS{\dot{\s}}}\\
&&\\
\qlb W_\Gb^\ci(\Lb,v) \rrto^{\qquad\DS{\Th}} && \AC\\
\enddiagram$$
$$\diagram
\qlb W_{\Gb_1}^\ci(\Lb_1,v_1) \rrto^{\qquad\DS{\Th^\pr_1}} \ddto_{\DS{\s}}&& \AC_1
\ddto^{\DS{\dot{\s}}}\\
&&\\
\qlb W_\Gb^\ci(\Lb,v) \rrto^{\qquad\DS{\Th'}} && \AC\\
\enddiagram\leqno{\mathrm{\it and}}$$
are commutative.
\end{prop}

\proof (a) follows from the Proper Base Change Theorem 
\cite[Chapter VI, Corollary 2.3]{Milne} applied 
to the Cartesian square 
$$\diagram
\Ybt_1 \rrto^{\DS{\pi_1}} \ddto_{\DS{\sigt}} && \Yb_1 \ddto^{\DS{\s}} \\
&&\\
\Ybt \rrto_{\DS{\pi}} &&\Yb.
\enddiagram$$
(b) follows from Lusztig's description of $\AC$. The commutativity of the first 
diagram in (c) follows from the fact that
$\s(\Cb_1^{\Gb_1})=\Cb^\Gb$ and from Theorem \ref{lusztig theta} 
while the commutativity of the second one follows from Proposition \ref{wouaw}.\fin

\bi

\begin{coro}\label{coro gamma}
We have $\g_{\Lb,v,\z}^\Gb\ci \s=\g_{\Lb_1,v_1,\z_1}^{\Gb_1}$.
\end{coro}

\bi

\sub{Parabolic restriction} Let 
$\Qb$ be a \para of $\Gb$ containing $\Pb$ and let 
$\Mb$ be the \levi of $\Qb$ containing $\Lb$. It follows from 
\cite[Theorem 8.3 (b)]{luicc} that~:

\begin{prop}\label{restriction parabolique}
$\g_{\Lb,v,\z}^\Mb=\Res_{W_\Mb^\ci(\Lb,v)}^{W_\Gb^\ci(\Lb,v)} \g_{\Lb,v,\z}^\Gb$.
\end{prop}

\bi

\remark{non restriction} 
If $\Gb'$ is a connected reductive subgroup of $\Gb$ which contains 
$\Lb$, then it may happen that 
$\g_{\Lb,v,\z}^{\Gb'}\not=
\Res_{W_{\Gb'}^\ci(\Lb,v)}^{W_\Gb^\ci(\Lb,v)} \g_{\Lb,v,\z}^\Gb$. 
An example is provided by the group $\Gb=\Sb\pb_4(\FM)$, as it will be shown in 
Part II of this paper.\finl

\bi



\rem Whenever $C_\Gb(u) \incl \Pb$, then Corollary \ref{coro zeta phi} shows that 
\ref{gamma produit}, Proposition \ref{K1}, and Proposition 
\ref{restriction parabolique} may be deduced from \ref{produit H}, 
Proposition \ref{G1 G H}, and Proposition \ref{prop compa para} respectively.\finl

\bi

\sec{Introducing Frobenius\label{part finite}}

\bi

\sub{Hypothesis and notation} 
In this section, and only in this section, we assume that $\FM$ is an 
algebraic closure of a finite field. In particular, $p > 0$. We fix a power 
$q$ of $p$ and we denote by $\fq$ the subfield of $\FM$ with $q$ elements. 
We assume also that $\Gb$ is defined over $\fq$ and we denote by 
$F : \Gb \to \Gb$ the corresponding Frobenius endomorphism. 
If $g \in \Gb^F$, we denote by $[g]$ (or $[g]_{\Gb^F}$ if 
necessary) the $\Gb^F$-conjugacy class of $g$.

We keep the notation introduced in \SEC\ref{sec endo} ($\Lb$, 
$\Cb$, $v$, $\EC$, $K$, $\Th$, $\g_{\Lb,v,\z}^\Gb$...). We assume that 
$\Lb$ is $F$-stable. Then, by Theorem \ref{rappellusztig} (e), there exists 
$n \in N_\Gb(\Lb)$ such that $F(\Pb)=\lexp{n}{\Pb}$. Now, by Lang 
theorem, we can pick an element $g \in \Gb$ such that $g^{-1}F(g)=n^{-1}$. 
Then $\lexp{g}{\Lb}$ and $\lexp{g}{\Pb}$ are $F$-stable. Since we are 
interested in the family of all $F$-stable subgroups of $\Gb$ which are 
conjugate to $\Lb$ under $\Gb$, we may, and we will assume that $\Lb$ and 
$\Pb$ are both $F$-stable. Without loss of generality, me may also 
assume that $\Bb$ and $\Tb$ are $F$-stable.

We also assume that $v$ and $\EC$ are $F$-stable.
Let $w \in W_\Gb^\ci(\Lb,v)$. We choose an element $g_w \in \Gb$ such 
that $g_w^{-1}F(g_w)=\dot{w}^{-1}$ (recall that $\dot{w}$ is a 
representative of $w$ in $N_\Gb(\Lb) \cap C_\Gb^\ci(v)$). We then put~:
$$\Lb_w=\lexp{g_w}{\Lb},\quad\quad v_w=\lexp{g_w}{v},
\quad\quad \Cb_w=\lexp{g_w}{\Cb},$$
$$\EC_w=(\ad g_w^{-1})^*\EC,\quad\quad{\mathrm{and}}
\quad\quad \FC_w=(\ad g_w^{-1})^*\FC.$$
Then $\Lb_w$ is an $F$-stable \levi of a \para of $\Gb$, $v_w \in \Lb_w^F$ 
is conjugate to $v$ in $\Gb^F$ (because $g_w^{-1}F(g_w) \in C_\Gb^\ci(v)$), 
$\Cb_w$ is the conjugacy class of $v_w$ in 
$\Lb_w$, $\EC_w$ is an $F$-stable cuspidal local system on $\Cb_w$ and 
$\FC_w=\qlb \boxtimes \EC_w$ (as a local system on 
$\Sigb_w=\Zb(\Lb_w)^\ci \times \Cb_w$).

\bi

\sub{A conjugacy result} In \cite[Proposition 2.1]{Bonnafe}, the author proved the 
following result~:

\bi

\begin{prop}\label{conjugue M}
Let $\Mb$ and $\Mb'$ be two $F$-stable \levis of $($non necessarily $F$-stable$)$ 
\paras of $\Gb$ which are geometrically conjugate and let $u$ be a unipotent element of 
$\Mb^F$. Assume the following conditions 
holds~:

\tete{a} $u$ is a distinguished element of $\Mb$,

\tete{b} $N_\Gb(\Mb)$ stabilizes the class $(u)_\Mb$, and

\tete{c} $A_\Mb(u) = A_\Gb(u)$.

\noindent Then $[u]_{\Gb^F} \cap \Mb'$ is a single $\Mb^{\pr F}$-conjugacy class.
\end{prop}

\bi

\begin{coro}\label{conjugue w}
If $w \in W_\Gb^\ci(\Lb,v)$, then $[v]_{\Gb^F} \cap \Lb_w^F=[v_w]_{\Lb_w^F}$.
\end{coro}

\bi

\proof This follows from Theorem \ref{rappellusztig} (a) and (c), from 
Theorem \ref{alv} and from the previous Proposition \ref{conjugue M}.\fin

\bi

\sub{Characteristic functions} We choose once and for all an isomorphism of local 
systems $\ph : F^*\EC \to \EC$ and we 
denote by $\XC_{\EC,\ph}$ the class function on $\Lb^F$ defined by
$$\XC_{\EC,\ph}(l)=\left\{\begin{array}{ll}
\Tr(\ph_l,\EC_l) & {\mathrm{if~}}l \in \Cb^F, \\
0 & {\mathrm{otherwise}},
\end{array}\right.$$
for any $l \in \Lb^F$. Using the isomorphism $\Th : \qlb W_\Gb^\ci(\Lb,v) \to \AC$, 
Lusztig defined  an isomorphism of local systems $\ph_w : F^*\EC_w \to \EC_w$. 
We recall his construction \cite[9.3]{lugf}. Let $\th_w : 
\EC \to (\ad \dot{w})^* \EC$ 
be the isomorphism of local systems defined in Theorem \ref{lusztig theta}. 
Then $\th_w$ induces an isomorphism of local systems 
$$F^* \ci (\ad g_w^{-1})^* \th_w : F^* \EC_w \longto (\ad g_w^{-1})^* \ci F^* \EC.$$
Moreover, $\ph$ induces an isomorphism 
$$(\ad g_w^{-1})^* \ph : (\ad g_w^{-1})^* \ci F^* \EC \longto \EC_w.$$
By composition of the two previous isomorphisms, we get an isomorphism
$$\ph_w : F^* \EC_w \longto \EC_w.$$
Once the isomorphism $\ph$ is chosen, the isomorphism $\ph_w$ depends only 
on the construction of the isomorphism $\th_w$. By Corollary 
\ref{epsilon}, the knowledge 
of the isomorphism $\th_w$ is equivalent to the knowledge of the linear character 
$\g_{\Lb,v,\z}^\Gb$.

\bi

%
%
If $\ch$ is an $F$-stable 
irreducible character of $W_\Gb^\ci(\Lb,v)$, then we denote by $\chit$ 
the {\it preferred extension} of $\ch$ to $W_\Gb^\ci(\Lb,v) \rtimes <F>$ 
(the preferred extension has been defined by Lusztig \cite{lucs}). 
The choice of $\ph$ and $\chit$ induces a well-defined isomorphism $\ph_\ch : 
F^*K_\ch \longmapright{\sim} K_\ch$. Then we set, for every $g \in \Gb^F$, 
$$\XC_{K_\ch,\ph_\ch}(g)=\left\{\begin{array}{ll}
\DS{\sum_{i \in \ZM} (-1)^i \Tr(\HC_g^i(\ph_\ch),\HC_g^i K_\ch)} & {\mathrm{if}}~g
~{\mathrm{is}~} unipotent,\\
0 & {\mathrm{otherwise}}
\end{array}\right.$$
The importance of the knowledge of the characteristic functions $\ch_{\EC_w,\ph_w}$ 
(therefore, of the linear character $\g_{\Lb,v,\z}^\Gb$) is given by 
the following theorem.

\bi

\begin{theo}[{\bf Lusztig}]\label{formule lusztig}
If $p$ is almost good for $\Gb$ and if $q$ is large enough, we have
$$\XC_{K_\ch,\ph_\ch} = {1 \over |W_\Gb^\ci(\Lb,v)|} 
\sum_{w \in W_\Gb^\ci(\Lb,v)} \chit(wF) R_{\Lb_w}^\Gb (\XC_{\EC_w,\ph_w}).$$
\end{theo}

\bi

\remark{large enough} The expression ``$q$ is large enough'' in Theorem 
\ref{formule lusztig} means that there exists a constant $q_0$ depending 
only on the root datum of $\Gb$ such that Theorem \ref{formule lusztig} 
holds if $q \ge q_0$. We will say that the pair $(\Gb,F)$ is {\it friendly} 
if the formula given by Theorem \ref{formule lusztig} holds in every 
connected reductive subgroup of $\Gb$ having the same rank.

Theorem \ref{formule lusztig} says that, if $p$ is almost good for $\Gb$ and 
if $q$ is large enough, then $(\Gb,F)$ is friendly.\finl 

\bi

\remark{independance} If $(\Gb,F)$ is friendly, then 
\cite[Proof of Theorem 6.1.1]{q grand} Mackey formula for Lusztig functors 
holds in $\Gb$. Consequently, the Lusztig functor 
$R_{\Lb_w \incl \Pb_w}^\Gb$ is known \cite[Proposition 7.1.1]{q grand} 
to be independent of the choice of 
the \para $\Pb_w$ of $\Gb$ which admits $\Lb_w$ as a Levi subgroup. 
That is why we denoted it simply by $R_{\Lb_w}^\Gb$.\finl

\bi

We conclude this section by determining explicitly the characteristic 
functions $\XC_{\EC_w,\ph_w}$ using the linear character $\g_{\Lb,v,\z}^\Gb$.
It follows from Lang's theorem that the set of rational conjugacy classes 
contained in $\Cb_w^F$ is in one-to-one correspondence with $H^1(F,A_\Lb(v_w))
\simeq H^1(\wdo F,A_\Lb(v)) =H^1(F,A_\Lb(v))$ (the last equality 
follows from the fact that $W_\Gb^\ci(\Lb,v)$ acts trivially 
on $A_\Lb(v)$). Let $a \in H^1(F,A_\Lb(v))$. We denote by $\aha$ 
a representative of $a$ in $A_\Lb(v)$ and by $v_{w,a}$ a representative 
of the rational conjugacy class contained in $\Cb_w^F$ parameterized by $a$. 
If $w=1$, we denote by $v_a$ the element $v_{w,a}$. We have $v_a \in \Lb^F$.  
It must be noticed that $[v_{w,a}]_{\Lb_w^F}=[v_a]_{\Gb^F} \cap \Lb^F$ 
(\cf Corollary \ref{conjugue w}).

By following step by step the construction of the \isos 
$\ph_w$, we obtain that the link between the class functions 
$\XC_{\EC,\ph}$ and $\XC_{\EC_w,\ph_w}$ 
is given in terms of the linear character $\g_{\Lb,v,\z}^\Gb$. 
More precisely, we get~:

\bi

\begin{prop}\label{fonction}
Let $w \in W_\Gb^\ci(\Lb,v)$ and let $a \in H^1(F,A_\Lb(v))$. Then 
$$\XC_{\EC_w,\ph_w}(v_{w,a})=\XC_{\EC,\ph}(v_a)\g_{\Lb,v,\z}^\Gb(w).$$
\end{prop}

\bi

Assume now until the end of this section that $\z$ is a linear character. 
In this case, we have
\equat\label{je sais pas}
\XC_{\EC_w,\ph_w}(v_{w,a})=\XC_{\EC_w,\ph_w}(v_w)\z(\aha).
\endequat
Note that $\z(\aha)$ does not depend on the choice of $\aha$ 
because $\z$ is $F$-stable. Hence, we deduce from Proposition 
\ref{fonction} the following result~:

\bi

\begin{coro}\label{description chi}
Assume that $\z$ is a linear character, and let $w \in W_\Gb^\ci(\Lb,v)$ 
and $a \in H^1(F,A_\Lb(v))$. Then 
$$\XC_{\EC_w,\ph_w}(v_{w,a})=
\XC_{\EC,\ph}(v)\g_{\Lb,v,\z}^\Gb(w)\z(\aha).$$
On the other hand, if $l \not\in \Cb_w^F$, then
$$\XC_{\EC_w,\ph_w}(l)=0.$$
\end{coro}

\bi

\remark{importance} In the theory of character sheaves applied to finite 
reductive groups, the characteristic functions $\XC_{\EC_w,\ph_w}$ play a 
crucial role, as it is shown in Theorem \ref{formule lusztig}. 
The Corollary \ref{description chi} shows the importance 
of the determination of the linear character $\g_{\Lb,v,\z}^\Gb$. 
We will show in Part II 
how the knowledge of the linear character $\g_{\Lb,v,\z}^\Gb$ 
whenever $v$ is regular and $p$ is good for $\Lb$ leads to an 
improvement of Digne, Lehrer and Michel's theorem on Lusztig 
restriction of Gel'fand-Graev characters \cite[Theorem 3.7]{DLM2}.\finl

\bi

\remark{scalaire} The characteristic function $\XC_{\EC_w,\ph_w}$ depends 
on the choice of the isomorphism $\ph$. Since $\EC$ is an irreducible 
local system, two \isos between $F^*\EC$ and $\EC$ differ only by 
a scalar. Hence the two characteristic functions they define 
differ also by the same constant. This shows that the formula 
in Corollary \ref{description chi} cannot be improved, because 
the factor $\XC_{\EC,\ph}(v)$ depends on the choice of the 
\iso $\ph$. We can give it, by multiplying $\ph$ by a scalar, 
any value we want.\finl

\bi

\bigskip

\bigskip

\hspace{2cm} C\'edric Bonnaf\'e

\hspace{2cm} CNRS - UMR 6623

\hspace{2cm} Universit\'e de Franche-Comt\'e

\hspace{2cm} D\'epartement de Math\'ematiques

\hspace{2cm} 16 Route de Gray

\hspace{2cm} 25000 BESAN\c{C}ON - FRANCE

\hspace{2cm} bonnafe\arobas math.univ-fcomte.fr

\end{document}